\documentclass{imanum}

\jno{drnxxx}
\received{\today}

\usepackage{graphicx}
\PassOptionsToPackage{table}{xcolor}

\newcommand{\eps}{\epsilon}

\renewcommand{\Re}{\mathbb{R}}           

\newcommand{\Iscr}{{\mathcal I}}

\newcommand{\beqn}[1]{\begin{equation}\label{#1}}
\newcommand{\eeqn}{\end{equation}}

\newcommand{\algbox}{%
  \vbox{\hrule height0.3pt\hbox{\vrule height2pt%
      width0.3pt\hskip\textwidth\hskip-0.6pt\vrule width0.3pt}\hrule height0.3pt}}
\newcommand{\AlgBegin}{\vspace{\smallskipamount}\noindent\algbox\par\vspace{\smallskipamount}\par\noindent}
\newcommand{\AlgEnd}{\par\noindent\algbox\vspace{\smallskipamount}}
\newcounter{pseudocode}[section]
\def\thepseudocode{\thesection.\arabic{pseudocode}}

\newcommand{\mystrut}{\vrule height9.5pt depth1.5pt width0pt}
\newcommand{\alg}[1]{\par\noindent\mystrut\ignorespaces\hbox to\textwidth{#1\hfill}}
\newcommand{\algt}[1]{\par\noindent\mystrut\hbox to\textwidth{\ignorespaces\hskip1.5em#1\hfill}}
\newcommand{\algtt}[1]{\par\noindent\mystrut\hbox to\textwidth{\ignorespaces\hskip3.5em#1\hfill}}
\newcommand{\algttt}[1]{\par\noindent\mystrut\hbox to\textwidth{\ignorespaces\hskip5.5em#1\hfill}}
\newcommand{\algtttt}[1]{\par\noindent\mystrut\hbox to\textwidth{\ignorespaces\hskip7.5em#1\hfill}}
\newcommand{\algttttt}[1]{\par\noindent\mystrut\hbox to\textwidth{\ignorespaces\hskip9.5em#1\hfill}}
\newcommand{\clg}[2]{\par\noindent\mystrut\hbox to\textwidth{\ignorespaces#1\hfill[#2]}}
\newcommand{\clgt}[2]{\par\noindent\mystrut\hbox to\textwidth{\ignorespaces\hskip1.5em#1\hfill[#2]}}
\newcommand{\clgtt}[2]{\par\noindent\mystrut\hbox to\textwidth{\ignorespaces\hskip3.5em#1\hfill[#2]}}
\newcommand{\clgttt}[2]{\par\noindent\mystrut\hbox to\textwidth{\ignorespaces\hskip5.5em#1\hfill[#2]}}
\newcommand{\clgtttt}[2]{\par\noindent\mystrut\hbox to\textwidth{\ignorespaces\hskip7.5em#1\hfill[#2]}}
\newcommand{\clgttttt}[2]{\par\noindent\mystrut\hbox to\textwidth{\ignorespaces\hskip9.5em#1\hfill[#2]}}

\usepackage{amsmath,amsfonts,multirow}
\newcommand{\Acal}{{\cal A}}

\newcommand{\Ical}{{\cal I}}

\newcommand{\Kcal}{{\cal K}}
\newcommand{\Lcal}{{\cal L}}

\newcommand{\Ocal}{{\cal O}}
\newcommand{\Pcal}{{\cal P}}

\newcommand{\Rcal}{{\cal R}}
\newcommand{\Scal}{{\cal S}}

\newcommand{\oline}[1]{\mkern 1.5mu\overline{\mkern-1.5mu#1}}

\renewcommand{\hbar}{\oline{h}}

\newcommand{\kbar}{\oline{k}}

\newcommand{\xbar}{\oline{x}}

\newcommand{\khat}{\hat{k}}

\usepackage{amssymb,algorithm}
\newtheorem{assumption}[theorem]{Assumption}
\newcommand{\balgorithm}  {\begin{algorithm}}
\newcommand{\ealgorithm}  {\end{algorithm}}
\newcommand{\balgorithmic}{\begin{algorithmic}}
\newcommand{\ealgorithmic}{\end{algorithmic}}
\newcommand{\bassumption} {\begin{assumption}}
\newcommand{\eassumption} {\end{assumption}}
\newcommand{\bcorollary}  {\begin{corollary}}
\newcommand{\ecorollary}  {\end{corollary}}
\newcommand{\bdefinition} {\begin{definition}}
\newcommand{\edefinition} {\end{definition}}
\newcommand{\blemma}      {\begin{lemma}}
\newcommand{\elemma}      {\end{lemma}}
\newcommand{\bproblem}    {\begin{problem}}
\newcommand{\eproblem}    {\end{problem}}
\newcommand{\bproof}      {\begin{proof}}
\newcommand{\eproof}      {\end{proof}}
\newcommand{\btheorem}    {\begin{theorem}}
\newcommand{\etheorem}    {\end{theorem}}

\DeclareMathOperator{\Null}     {Null}

\DeclareMathOperator{\st}       {s.\!t.}

\newcommand{\thalf}{\tfrac12}
\renewcommand{\(}{\left(}
\renewcommand{\)}{\right)}

\newcommand{\N}[1]{\mathbb{N}^{#1}}

\newcommand{\R}[1]{\mathbb{R}^{#1}}

\newcommand{\baligned}     {\begin{aligned}}
\newcommand{\ealigned}     {\end{aligned}}
\newcommand{\barray}       {\begin{array}}
\newcommand{\earray}       {\end{array}}
\newcommand{\bbmatrix}     {\begin{bmatrix}}
\newcommand{\ebmatrix}     {\end{bmatrix}}
\newcommand{\bcases}       {\begin{cases}}
\newcommand{\ecases}       {\end{cases}}
\newcommand{\bcenter}      {\begin{center}}
\newcommand{\ecenter}      {\end{center}}
\newcommand{\bcolumn}      {\begin{column}}
\newcommand{\ecolumn}      {\end{column}}
\newcommand{\bcolumns}     {\begin{columns}}
\newcommand{\ecolumns}     {\end{columns}}
\newcommand{\benumerate}   {\begin{enumerate}}
\newcommand{\eenumerate}   {\end{enumerate}}
\newcommand{\bequation}    {\begin{equation}}
\newcommand{\eequation}    {\end{equation}}
\newcommand{\bequationn}   {\begin{equation*}}
\newcommand{\eequationn}   {\end{equation*}}
\newcommand{\bfigure}      {\begin{figure}}
\newcommand{\efigure}      {\end{figure}}
\newcommand{\bflushright}  {\begin{flushright}}
\newcommand{\eflushright}  {\end{flushright}}
\newcommand{\bitemize}     {\begin{itemize}}
\newcommand{\eitemize}     {\end{itemize}}
\newcommand{\bpmatrix}     {\begin{pmatrix}}
\newcommand{\epmatrix}     {\end{pmatrix}}
\newcommand{\bsubequations}{\begin{subequations}}
\newcommand{\esubequations}{\end{subequations}}
\newcommand{\btable}       {\begin{table}}
\newcommand{\etable}       {\end{table}}
\newcommand{\btabular}     {\begin{tabular}}
\newcommand{\etabular}     {\end{tabular}}
\newcommand{\bvmatrix}     {\begin{vmatrix}}
\newcommand{\evmatrix}     {\end{vmatrix}}
\newcommand{\bequalin}     {\bequationn\baligned}
\newcommand{\eequalin}     {\ealigned\eequationn}
\newcommand{\bequali}      {\bsubequations\begin{align}}
\newcommand{\eequali}      {\end{align}\esubequations}

\usepackage[noend]{algpseudocode}
\usepackage{xcolor}
\usepackage{longtable}
\usepackage{multirow}
\usepackage{multicol}
\usepackage{afterpage}

\newcommand{\ARC}        {\textsc{arc}}
\newcommand{\iARC}       {\texttt{iARC}}
\newcommand{\iRNewton}   {\texttt{iR\_Newton}}
\newcommand{\TRACE}      {\textsc{trace}}
\newcommand{\gammal}     {\gamma_1}
\newcommand{\gammau}     {\gamma_2}
\newcommand{\sigmamin}   {\underline\sigma}
\newcommand{\sigmamax}   {\overline\sigma}
\newcommand{\sigmal}     {\sigma^{\textsc{l}}}
\newcommand{\sigmau}     {\sigma^{\textsc{u}}}

\newcommand{\betal}      {\beta^{\textsc{l}}}
\newcommand{\betau}      {\beta^{\textsc{u}}}
\newcommand{\betan}      {\beta^{\textsc{n}}}

\algrenewcommand{\algorithmiccomment}[1]{\hfill[#1]}
\algnewcommand{\LineComment}{\Statex \hskip\ALG@thistlm}
\newcommand\AlgBreak{\vspace*{-.7\baselineskip}\Statex\hspace*{\dimexpr-\algorithmicindent-2pt\relax}\rule{\textwidth}{0.4pt}}

\begin{document}

\title {An Inexact Regularized Newton Framework with a Worst-Case Iteration Complexity of $\Ocal(\epsilon^{-3/2})$ for Nonconvex Optimization\thanks{This work was supported in part by the U.S.~Department of Energy, Office of Science, Applied Mathematics, Early Career Research Program under Award Number DE--SC0010615 and by the U.S.~National Science Foundation, Division of Mathematical Sciences, Computational Mathematics Program under Award Number DMS--1016291.}}
\shorttitle{An Inexact Regularized Newton Framework with Complexity $\Ocal(\epsilon^{-3/2})$}
\author{{\sc Frank~E.~Curtis}\thanks{Corresponding author. Email: frank.e.curtis@gmail.com}\\[2pt]
        Department of Industrial and Systems Engineering, Lehigh University\\[6pt]
        {\sc and}\\[6pt]
        {\sc Daniel~P.~Robinson}\thanks{Email: daniel.p.robinson@gmail.com}\\[2pt]
        Department of Applied Mathematics and Statistics, Johns Hopkins University\\[6pt]
        {\sc and}\\[6pt]
        {\sc Mohammadreza Samadi}\thanks{Email: mohammadreza.samadi@lehigh.edu}\\[2pt]
        Department of Industrial and Systems Engineering, Lehigh University
       }
\shortauthorlist{F.~E.~Curtis \emph{et al.}}

\maketitle

\begin{abstract}
  {An algorithm for solving smooth nonconvex optimization problems is proposed that, in the worst-case, takes $\Ocal(\epsilon^{-3/2})$ iterations to drive the norm of the gradient of the objective function below a prescribed positive real number $\epsilon$ and can take $\Ocal(\eps^{-3})$ iterations to drive the leftmost eigenvalue of the Hessian of the objective above $-\eps$.  The proposed algorithm is a general framework that covers a wide range of techniques including quadratically and cubically regularized Newton methods, such as the Adaptive Regularisation using Cubics (\ARC) method and the recently proposed Trust-Region Algorithm with Contractions and Expansions (\TRACE).  The generality of our method is achieved through the introduction of generic conditions that each trial step is required to satisfy, which in particular allow for inexact regularized Newton steps to be used.  These conditions center around a new subproblem that can be approximately solved to obtain trial steps that satisfy the conditions.  A new instance of the framework, distinct from \ARC{} and \TRACE{}, is described that may be viewed as a hybrid between quadratically and cubically regularized Newton methods.  Numerical results demonstrate that our hybrid algorithm outperforms a cublicly regularized Newton method.}
{unconstrained optimization, nonlinear optimization, nonconvex optimization, inexact Newton methods, worst-case iteration complexity, worst-case evaluation complexity}
\end{abstract}

\section{Introduction}\label{sec.introduction}

This paper proposes an algorithm for solving
\bequation \label{prob.f}
  \min_{x \in \R{n}}\ f(x),
\eequation
where the (possibly nonconvex) objective function $f:\R{n} \to \R{}$ is assumed to be twice-continuously differentiable. The optimization problem~\eqref{prob.f} has been widely studied, as evidenced by its appearance as the focal point of numerous textbooks; e.g., see~\cite{BazaSherShet06}, \cite{Bert99}, \cite{ConnGoulToin00}, \cite{GrivNashSofe08}, \cite{NoceWrig06}, and \cite{Rusz06}. 

For many years, the most popular methods for solving~\eqref{prob.f} were in classes known as line search and trust region methods.  Recently, however, cubic regularization methods have become popular, which are based on the pioneering work by \cite{Grie81} and \cite{NestPoly06}. Their rise in popularity is due to increased interest in algorithms with improved complexity properties, which stems from the impact of so-called optimal algorithms for solving convex optimization problems.  For problem~\eqref{prob.f}, by complexity properties, we mean a guaranteed bound on the number of iterations (or function evaluations or derivative evaluations) needed by an algorithm before the norm of the gradient of the objective must fall below a positive threshold $\epsilon > 0$.  In other words, if $x_k$ denotes the $k$th iteration of an algorithm, one seeks a bound on the number of iterations until it is guaranteed that
\bequationn
  \|\nabla f(x_k)\| \leq \epsilon.
\eequationn
The complexity of a traditional trust region method~(e.g., see Algorithm 6.1.1 in \cite{ConnGoulToin00}) is $\Ocal(\epsilon^{-2})$ (see \cite{CartGoulToin10}), which falls short of the $\Ocal(\epsilon^{-3/2})$ complexity for cubic regularization methods (e.g., see the \ARC{} method by~\cite{CartGoulToin11a,CartGoulToin11b}).  This latter complexity is optimal among a certain broad class of second-order methods when employed to minimize a broad class of objective functions; see \cite{CartGoulToin11c}.  That said, one can obtain even better complexity properties if higher-order derivatives are used; see~\cite{BirgGardMartSantToin17} and \cite{CartGouToin17}.

The better complexity properties of regularization methods such as \ARC\ have been a major point of motivation for discovering other methods that attain the same worst-case iteration complexity bounds.  For example, the recently introduced (nontraditional) trust region method known as \TRACE\ (see~\cite{CurtRobiSama17}) has the same optimal $\Ocal(\epsilon^{-3/2})$ complexity, while at the same time allowing traditional trust region trial steps to be computed and used.  A key aspect of the \TRACE{} framework is that a solution to an implicit trust region problem is obtained by varying a regularization parameter instead of a trust region radius. This key idea has been adopted and advanced further by~\cite{BirgMart17_siopt}; in particular, they propose an algorithm that has optimal iteration complexity by solving quadratic subproblems that have a carefully chosen quadratic regularization parameter.

\paragraph{Contributions}  The main contributions of this paper relate to advancing the understanding of optimal complexity algorithms for solving the smooth optimization problem~\eqref{prob.f}. Our proposed framework is intentionally very general; it is not a trust region method, a quadratic regularization method, or a cubic regularization method.  Rather, we propose a generic set of conditions that each trial step must satisfy that still allow us to establish an optimal first-order complexity result as well as a second-order complexity bound similar to the methods above.  Our framework contains as special cases other optimal complexity algorithms such as \ARC{} and \TRACE{}.  To highlight this generality of our contribution, we describe one particular instance of our framework that appears to be new to the literature.

During the final preparation of this article, we came across the work in~\cite{Duss17} and \cite{DussOrba17}.  This work shares certain commonalities with our own and appears to have been developed at the same time.  Although there are numerous differences, we shall only point out three of them. First, the precise conditions that they require for each trial step are different from ours.  In particular, the condition stated as (3.1c) in~\cite{DussOrba17} requires that regularization is used to compute every trial step, a property not shared by our method (which can employ Newton steps).  Second, they do not consider second-order convergence or complexity properties, although they might be able to do so by incorporating second-order conditions similar to ours.  Third, they focus on strategies for identifying an appropriate value for the regularization parameter.  An implementation of our method might consider their proposals, but could employ other strategies as well.  In any case, overall, we believe that our papers are quite distinct, and in some ways are complementary.

\paragraph{Organization}  In~\S\ref{sec.algorithm}, we present our general framework that is formally stated as Algorithm~\ref{alg.inexact_regularized_Newton}. In~\S\ref{sec.convergence}, we prove that our framework enjoys first-order convergence (see \S\ref{sec.global}), an optimal first-order complexity (see \S\ref{sec.complexity}), and certain second-order convergence and complexity guarantees (see~\S\ref{sec.second_order}).  In~\S\ref{sec.instances}, we show that \ARC{} and \TRACE{} can be viewed as special cases of our framework, and present yet another instance that is distinct from these methods.  In~\S\ref{sec.numerical}, we present details of implementations of a cubic regularization method and our newly proposed instance of our framework, and provide the results of numerical experiments with both.  Finally, in~\S\ref{sec.conclusion}, we present final comments.

\paragraph{Notation}  We use $\R{}_+$ to denote the set of nonnegative scalars, $\R{}_{++}$ to denote the set of positive scalars, and $\N{}_+$ to denote the set of nonnegative integers.  Given a real symmetric matrix $A$, we write $A \succeq 0$ (respectively, $A \succ 0$) to indicate that $A$ is positive semidefinite (respectively, positive definite).  Given a pair of scalars $(a,b) \in \R{} \times \R{}$, we write $a \perp b$ to indicate that $ab=0$.  Similarly, given such a pair, we denote their maximum as $\max\{a,b\}$ and their minimum as $\min\{a,b\}$.  Given a vector $v$, we denote its (Euclidean) $\ell_2$-norm as $\|v\|$. Finally, given a discrete set $\Scal$, we denote its cardinality by $|\Scal|$.

Corresponding to the objective $f : \R{n} \to \R{}$, we define the gradient function $g := \nabla f : \R{n} \to \R{n}$ and the Hessian function $H := \nabla^2 f : \R{n} \to \R{n \times n}$.  Given an iterate $x_k$ in an algorithm for solving \eqref{prob.f}, we define $f_k := f(x_k)$, $g_k := g(x_k) := \nabla f(x_k)$, and $H_k := H(x_k) := \nabla^2 f(x_k)$.  Similarly, we apply a subscript to other algorithmic quantities whose definition depends on the iteration number~$k$.

\section{Algorithm Description}\label{sec.algorithm}

Our algorithm involves generic conditions that a trial step toward solving problem~\eqref{prob.f} must satisfy.  One can obtain a step satisfying these conditions by computing---for appropriate positive lower and upper bounds $\sigmal_k$ and $\sigmau_k$, respectively, on the ratio between a regularization variable $\lambda \geq 0$ and the norm of the trial step---an approximate solution of the subproblem
\bequation\label{prob.primal-dual}
  \baligned
    \Pcal_k(\sigmal_k,\sigmau_k)\ :\ \min_{(s,\lambda) \in \R{n}\times\R{}_{+}} &\ f_k + g_k^T s + \thalf s^T (H_k + \lambda I) s \\
    \st &\ (\sigmal_k)^2\|s\|^2 \leq \lambda^2 \leq (\sigmau_k)^2\|s\|^2.
  \ealigned
\eequation
For a given value of the regularization variable $\lambda$, this problem involves a quadratic objective function and an upper bound on the norm of the trial step, just as in a trust region method.  However, it also includes a lower bound on the norm of the trial step, and, in general, with $\lambda$ as a variable, it encapsulates other types of subproblems as well, including those present in a cubic regularization framework.  For additional details on the properties of this subproblem and its solutions, see Appendices~\ref{app.subproblem} and \ref{app.subproblem.reduced}.

The conditions that the $k$th trial step and regularization pair, i.e., $(s_k,\lambda_k)$, must satisfy are stated in Assumption~\ref{ass.subproblem_accuracy} below, wherein we invoke the following (unregularized) quadratic model of $f$ at $x_k$:
\bequationn
  q_k(s) := f_k + g_k^T s + \thalf s^T H_k s.
\eequationn

\bassumption\label{ass.subproblem_accuracy}
  \textit{
  The pair $(s_k,\lambda_k)$ is computed such that it is feasible for problem~\eqref{prob.primal-dual} and, with
  \bequation\label{def.delta}
    \Delta_k(s_k,\lambda_k) := \bcases \|s_k\| & \text{if $\lambda_k = 0$} \\ \frac{1}{\sqrt{6}}\sqrt{\frac{\|g_k\|\|s_k\|}{\lambda_k}} & \text{if $\lambda_k > 0$} \ecases
  \eequation
  and constants $(\kappa_1,\kappa_2,\kappa_3) \in \R{}_{++} \times \R{}_{++} \times \R{}_{++}$, the following hold:
  \bsubequations\label{eq.subproblem_accuracy}
    \begin{align}
      f_k - q_k(s_k)
        &\geq \frac{\|g_k\|}{6\sqrt{2}} \min \left\{ \frac{\|g_k\|}{1+\|H_k\|}, \Delta_k(s_k,\lambda_k) \right\}; \label{eq.Cauchy_decrease} \\
      s_k^T (g_k + (H_k+\lambda_k I)s_k)
        &\leq \min\{\kappa_1 \|s_k\|^2,\thalf s_k^T (H_k+\lambda_k I)s_k + \thalf \kappa_2 \|s_k\|^3\};\ \ \text{and} \label{eq.subspace_optimality} \\
      \| g_k + (H_k + \lambda_k I) s_k \| &\leq \lambda_k \|s_k\| + \kappa_3 \|s_k\|^2. \label{eq.residual_error}
  \end{align}
  \end{subequations}
  }
\eassumption

To see that Assumption~\ref{ass.subproblem_accuracy} is well-posed and consistent with problem~\eqref{prob.primal-dual}, we refer the reader to Theorem~\ref{th.consistency-subspace} in Appendix~\ref{app.subproblem.reduced} wherein we prove that any solution of problem~\eqref{prob.primal-dual} with $s$ restricted to a sufficiently large dimensional subspace of $\R{n}$ satisfies all of the conditions in Assumption~\ref{ass.subproblem_accuracy}.  We also claim that one can obtain a pair satisfying Assumption~\ref{ass.subproblem_accuracy} in either of the following two ways:
\bitemize
  \item Choose $\sigma \in [\sigmal_k,\sigmau_k]$, compute $s_k$ by minimizing the cubic function
  \bequation\label{eq.c}
    c_k(s;\sigma) := q_k(s) + \thalf\sigma\|s\|^3 = f_k + g_k^Ts + \thalf s^TH_ks + \thalf\sigma\|s\|^3
  \eequation
  over a sufficiently large dimensional subspace of $\R{n}$ (assuming, when $\sigma = \sigmal_k = 0$, that this function is not unbounded below), then set $\lambda_k \gets \sigma\|s_k\|$.  This is essentially the strategy employed in cubic regularization methods such as \ARC{}.
  \item Choose $\lambda_k \geq 0$, then compute $s_k$ by minimizing the objective of~\eqref{prob.primal-dual} with $\lambda = \lambda_k$ over a sufficiently large dimensional subspace of $\R{n}$ (assuming that the function is not unbounded below).  The resulting pair $(s_k,\lambda_k)$ satisfies Assumption~\ref{ass.subproblem_accuracy} as long as it is feasible for \eqref{prob.primal-dual}.  This is essentially the strategy employed in \cite{BirgMart17_siopt} and partly employed in \TRACE{}.
\eitemize
One can imagine other approaches as well.  Overall, we state problem~\eqref{prob.primal-dual} as a guide for various techniques for computing the pair $(s_k,\lambda_k)$.  Our theory simply relies on the fact that any such computed pair satisfies the conditions in Assumption~\ref{ass.subproblem_accuracy}.

Our algorithm, stated as Algorithm~\ref{alg.inexact_regularized_Newton}, employs the following ratio (also employed, e.g., in \TRACE{}) to determine whether a given trial step is accepted or rejected:
\bequationn
  \rho_k := \frac{f_k - f(x_k + s_k)}{\|s_k\|^3}.
\eequationn
One potential drawback of employing this ratio is that the ratio is not invariant to scaling of the objective function.  However, the use of this ratio can still be justified.  For example, if one were to compute $s_k$ by minimizing the cubic model \eqref{eq.c} for some $\sigma > 0$, then the reduction in this model yielded by $s_k$ is bounded below by a fraction of $\sigma\|s_k\|^3$ (see \cite[Lemma~4.2]{CartGoulToin11b}), meaning that $\rho_k \geq \eta$ holds when $\sigma \geq \eta$ and the actual reduction in $f$ is proportional to the reduction in the cubic model.  For further justification for this choice---such as how it allows the algorithm to accept Newton steps when the norm of the trial step is small (and, indeed, the norms of accepted steps vanish asymptotically as shown in Lemma~\ref{lem.infinite_accepted} later on)---we refer the reader to \cite{BirgMart17_siopt} and \cite{CurtRobiSama17}.

\balgorithm[ht]
  \caption{Inexact Regularized Newton Framework}
  \label{alg.inexact_regularized_Newton}
  \balgorithmic[1]
    \smallskip
    \Require an acceptance constant $\eta \in \R{}_{++}$ with $0 < \eta < 1$
    \Require bound update constants $\{\gammal,\gammau\} \subset \R{}_{++}$ with $1 < \gammal \leq \gammau$
    \Require ratio lower and upper bound constants $\{\sigmamin,\sigmamax\} \subset \R{}_{++}$ such that $\sigmamax \geq \sigmamin$
    \smallskip
    \AlgBreak
    \Procedure{Inexact Regularized Newton}{}
      \State set $x_0 \in \R{n}$
      \State set $\sigmal_0 \gets 0$ and $\sigmau_0 \in [\sigmamin,\sigmamax]$ \label{step.initialize} 
      \For{$k \in \N{}_+$}
        \State set $(s_k,\lambda_k)$ satisfying Assumption~\ref{ass.subproblem_accuracy} \label{step.Pk}
        \If{$\rho_k \geq \eta$} \label{step.accept} \Comment{accept step}
          \State set $x_{k+1} \gets x_k + s_k$
          \State set $\sigmal_{k+1} \gets 0$ and $\sigmau_{k+1} \gets \sigmau_k$ \label{step.accept.next}
        \Else\ (i.e., $\rho_k < \eta$) \label{step.reject} \Comment{reject step}
          \State set $x_{k+1} \gets x_k$
          \If{$\lambda_k < \sigmamin \|s_k\|$} \label{step.reject.lambda} 
            \State set $\sigmal_{k+1} \in [\sigmamin,\sigmamax]$ and $\sigmau_{k+1} \in [\sigmal_{k+1},\sigmamax]$ \label{step.reject.lambda_small}
          \Else
            \State set $\sigmal_{k+1} \gets \gammal \frac{\lambda_k}{\|s_k\|}$ and $\sigmau_{k+1} \gets \gammau \frac{\lambda_k}{\|s_k\|}$ \label{step.reject.lambda_large}
          \EndIf
        \EndIf
      \EndFor
    \EndProcedure
    \smallskip
  \ealgorithmic
\ealgorithm

\section{Convergence Analysis}\label{sec.convergence}

In this section, we prove global convergence guarantees for Algorithm~\ref{alg.inexact_regularized_Newton}.  In particular, we prove under common assumptions that, from remote starting points, the algorithm converges to first-order stationarity, has a worst-case iteration complexity to approximate first-order stationarity that is on par with the methods in \cite{CartGoulToin11b}, \cite{CurtRobiSama17}, and \cite{BirgMart17_siopt}, and---at least in a subspace determined by the search path of the algorithm---converges to second-order stationarity with a complexity on par with the methods in \cite{CartGoulToin11b} and \cite{CurtRobiSama17}.

\subsection{First-Order Global Convergence}\label{sec.global}

Our goal in this subsection is to prove that the sequence of objective gradients vanishes.  We make the following assumption about the objective function, which is assumed to hold throughout this section.

\bassumption\label{ass.f}
  \textit{
  The objective function $f : \R{n} \to \R{}$ is twice continuously differentiable and bounded below by a scalar $f_{\inf} \in \R{}$ on $\R{n}$.
  }
\eassumption

We also make the following assumption related to the sequence of iterates.

\bassumption\label{ass.global}
  \textit{
  The gradient function $g : \R{n} \to \R{n}$ is Lipschitz continuous with Lipschitz constant $g_{Lip} \in \R{}_{++}$ in an open convex set containing the sequences $\{x_k\}$ and $\{x_k + s_k\}$.  Furthermore, the gradient sequence $\{g_k\}$ has $g_k \neq 0$ for all $k \in \N{}_+$ and is bounded in that there exists a scalar constant $g_{max}\in \R{}_{++}$ such that $\|g_k\| \leq g_{max}$ for all $k \in \N{}_+$.
  }
\eassumption

It is worthwhile to note in passing that our complexity bounds for first- and second-order stationarity remain true even if one were to consider the possibility that $g_k = 0$ for some $k \in \N{}_+$, in which case one would have the algorithm terminate finitely or, if $H_k \not\succeq 0$, compute an improving direction of negative curvature for $H_k$.  However, allowing this possibility---which is typically unlikely ever to occur in practice---would only serve to obscure certain aspects of our analysis.  We refer the reader, e.g., to \cite{CartGoulToin11b} (specifically, to the discussions at the ends of \S2.1, \S4, and \S5 in that work) for commentary about why zero gradient values do not ruin complexity guarantees such as we present.

We begin with two lemmas each revealing an important consequence of Assumptions~\ref{ass.f} and~\ref{ass.global}.

\blemma\label{lem.s_nonzero}
  \textit{
  For all $k \in \N{}_+$, it follows that $s_k \neq 0$.
  }
\elemma
\bproof
  The result follows by combining that $g_k \neq 0$ for all $k \in \N{}_+$ (see Assumption~\ref{ass.global}) with \eqref{eq.residual_error}.
\eproof

\blemma\label{lem.H_bounded}
  \textit{
  The Hessian sequence $\{H_k\}$ is bounded in norm in that there exists a scalar constant $H_{max}\in \R{}_{++}$ such that $\|H_k\| \leq H_{max}$ for all $k \in \N{}_+$.
  }
\elemma
\bproof
  The result follows by Assumption~\ref{ass.f}, the Lipschitz continuity of $g$ in Assumption~\ref{ass.global}, and Lemma~1.2.2 in \cite{Nest04}.
\eproof

In our next lemma, we prove an upper bound for the regularization variable $\lambda_k$.

\blemma\label{lem.lambda-upperbound}
  \textit{
  For all $k \in \N{}_+$, the pair $(s_k,\lambda_k)$ satisfies
  \bequationn
    \lambda_k \leq 2 \frac{\|g_k\|}{\|s_k\|} + \tfrac32 H_{max} + \kappa_1.
  \eequationn
  }
\elemma
\bproof
  Since \eqref{eq.Cauchy_decrease} ensures $q_k(s_k) - f_k \leq 0$, it follows with \eqref{eq.subspace_optimality} and Lemma~\ref{lem.H_bounded} that
  \bequalin
    0 \geq q_k(s_k) - f_k &= g_k^T s_k + \thalf s_k^T H_k s_k \\
      &\geq g_k^T s_k + \thalf s_k^T H_k s_k + s_k^T(g_k + (H_k+\lambda_k I) s_k) - \kappa_1\|s_k\|^2 \\
      &= 2g_k^T s_k + \tfrac32 s_k^T H_k s_k + \lambda_k \|s_k\|^2 - \kappa_1\|s_k\|^2 \\
      &\geq -2 \|g_k\| \|s_k\| - \tfrac32 H_{max}\|s_k\|^2 + \lambda_k \|s_k\|^2 - \kappa_1\|s_k\|^2.
  \eequalin
  After rearrangement and dividing by $\|s_k\|^2 \neq 0$ (see Lemma~\ref{lem.s_nonzero}), the desired result follows.
\eproof

Using Lemma~\ref{lem.lambda-upperbound}, we now prove a lower bound for the reduction in $q_k$ yielded by $s_k$.

\blemma\label{lem.Cauchydecrease}
  \textit{
  For all $k \in \N{}_+$, the step $s_k$ satisfies
  \bequationn
    f_k - q_k(s_k) \geq \frac{\|g_k\|}{6\sqrt{2}} \min \left\{\frac{\|g_k\|}{1+H_{max}},\frac{\|s_k\|}{\sqrt{6}}\sqrt{\frac{\|g_k\|}{2\|g_k\|+ \|s_k\|(\tfrac32 H_{max} + \kappa_1)}}\right\}.
  \eequationn
  }
\elemma
\bproof
  If $\lambda_k = 0$, then by \eqref{eq.Cauchy_decrease} and Lemma~\ref{lem.H_bounded} it follows that
  \bequationn
    f_k - q_k(s_k) \geq \frac{\|g_k\|}{6\sqrt{2}} \min \left\{\frac{\|g_k\|}{1+\|H_k\|}, \|s_k\| \right\} \geq \frac{\|g_k\|}{6\sqrt{2}} \min \left\{\frac{\|g_k\|}{1+H_{max}}, \|s_k\| \right\}.
  \eequationn
  On the other hand, if $\lambda_k > 0$, then \eqref{eq.Cauchy_decrease}, Lemma~\ref{lem.H_bounded}, and Lemma~\ref{lem.lambda-upperbound} imply that
  \bequalin 
    f_k - q_k(s_k)
      &\geq \frac{\|g_k\|}{6\sqrt{2}} \min \left\{\frac{\|g_k\|}{1+\|H_k\|}, \frac{1}{\sqrt{6}}\sqrt{\frac{\|g_k\|\|s_k\|}{\lambda_k}} \right\} \\
      &\geq \frac{\|g_k\|}{6\sqrt{2}} \min \left\{\frac{\|g_k\|}{1+H_{max}}, \frac{\|s_k\|}{\sqrt{6}}\sqrt{\frac{\|g_k\|}{2\|g_k\|+ \|s_k\|(\tfrac32 H_{max} + \kappa_1)}} \right\}.
  \eequalin
  Combining the inequalities from these two cases, the desired result follows.
\eproof

Going forward, for ease of reference, we respectively define sets of indices corresponding to accepted and rejected steps throughout a run of the algorithm as
\bequationn
  \Acal := \{k \in \N{}_+ : \rho_k \geq \eta\} \ \ \text{and} \ \ \Rcal := \{k \in \N{}_+ : \rho_k < \eta\}.
\eequationn
We now show that if the algorithm were only to compute rejected steps from some iteration onward, then the sequence $\{\lambda_k / \|s_k\|\}$ diverges to infinity.

\blemma\label{lem.ratio_diverges}
  \textit{
  If $k\in\Rcal$ for all sufficiently large $k \in \N{}_+$, then $\{\lambda_k / \|s_k\|\} \to \infty$.
  }
\elemma
\bproof
  Without loss of generality, assume that $\Rcal = \N{}_+$.  We now prove that the condition in Step~\ref{step.reject.lambda} cannot be true more than once. Suppose, in iteration $\khat \in \N{}_+$, Step~\ref{step.reject.lambda_small} is reached, which means that $\lambda_{\khat+1} / \|s_{\khat+1}\| \geq \sigmamin$ since $(s_{\khat+1},\lambda_{\khat+1})$ is required to be feasible for $\Pcal_{\khat+1}(\sigmal_{k+1},\sigmau_{k+1})$ in Step~\ref{step.Pk} where $\sigmal_{k+1} \geq \sigmamin$. Therefore, the condition in Step~\ref{step.reject.lambda} tests false in iteration $(\khat+1)$. Then, from Step~\ref{step.Pk}, Step~\ref{step.reject.lambda_large}, and the fact that $\gamma_1 > 1$, it follows that $\{\lambda_k / \|s_{k}\|\}$ is monotonically increasing for all $k \geq \khat$.  Therefore, the condition in Step~\ref{step.reject.lambda} cannot test true for any $k \geq \khat+1$.  Now, to see that the sequence diverges, notice from this fact, Step~\ref{step.Pk}, and Step~\ref{step.reject.lambda_large}, it follows that for all $k \geq \khat+1$ we have $\lambda_{k+1} / \|s_{k+1}\| \geq \gammal (\lambda_k / \|s_k\|)$ where $\gammal > 1$.  Thus, $\{\lambda_k / \|s_{k}\|\} \to \infty$, as claimed.
\eproof

We now prove that if the gradients are bounded away from zero and the sequence of ratios $\{\lambda_k / \|s_k\|\}$ diverges, then $\rho_k \geq \eta$ for all sufficiently large $k \in \N{}_+$, meaning that the steps are accepted.

\blemma \label{lem.eventuallyaccepted}
  \textit{
  Suppose that $\Ical \subseteq \N{}_+$ is an infinite index set such that for $\epsilon \in \R{}_{++}$ independent of $k$, one finds that $\|g_k\| \geq \epsilon$ for all $k \in \Ical$ and $\{\lambda_k/\|s_k\|\}_{k \in \Ical} \to \infty$. Then, for all sufficiently large $k \in \Ical$, it follows that $\rho_k \geq \eta$, meaning $k \in \Acal$.
  }
\elemma
\bproof
  From the Mean Value Theorem, there exists $\xbar_k \in [x_k,x_k+s_k]$ such that
  \begin{align}
    q_k(s_k) - f(x_k+s_k) 
      &= (g_k - g(\xbar_k))^Ts_k + \thalf s_k^T H_k s_k \nonumber \\
      &\geq -\|g_k - g(\xbar_k)\|\|s_k\| - \thalf \|H_k\|\|s_k\|^2. \label{eq.mvt}
  \end{align}
  From this, Lemma~\ref{lem.Cauchydecrease}, and Assumption~\ref{ass.global}, it follows that, for all $k \in \Ical$,
  \bequalin
    f_k - f(x_k + s_k)
       =&\ f_k - q_k(s_k) + q_k(s_k) - f(x_k+s_k)\\
    \geq&\ \frac{\|g_k\|}{6\sqrt{2}} \min \left\{\frac{\|g_k\|}{1+H_{max}},\frac{\|s_k\|  }{\sqrt{6}}\sqrt{\frac{\|g_k\|}{2\|g_k\|+ \|s_k\|(\tfrac32 H_{max} + \kappa_1)}}\right\} - (g_{Lip} + \thalf H_{max}) \|s_k\|^2 \\
    \geq&\ \frac{\epsilon}{6\sqrt{2}} \min \left\{\frac{\epsilon}{1+H_{max}},\frac{\|s_k\|}{\sqrt{6}}\sqrt{\frac{\epsilon}{2g_{max}+ \|s_k\|(\tfrac32 H_{max} + \kappa_1)}}\right\} - (g_{Lip} + \thalf H_{max}) \|s_k\|^2.
  \eequalin
  This shows that there exists a threshold $s_{thresh} > 0$ such that 
  \bequationn
    f_k - f(x_k+s_k) \geq \eta \|s_k\|^3\ \ \text{whenever}\ \ k \in \Ical\ \ \text{and}\ \ \|s_k\| \leq s_{thresh}.
  \eequationn 
  We now claim that $\{\|s_k\|\}_{k \in \Ical} \to 0$.  To prove this claim, suppose by contradiction that there exists an infinite subsequence $\Ical_s \subseteq\Ical$ and scalar $\epsilon_s \in \R{}_{++}$ such that $\|s_k\| \geq \epsilon_s$ for all $k \in \Ical_s$.  It then follows from the boundedness of $\{\|g_k\|\}$ (see Assumption~\ref{ass.global}) and Lemma~\ref{lem.lambda-upperbound} that $\{\lambda_k\}_{k\in\Ical_s}$ is bounded.  This allows us to conclude that $\{\lambda_k/\|s_k\|\}_{k \in \Ical_s}$ is bounded, which contradicts the assumptions of the lemma.  Thus, $\{\|s_k\|\}_{k \in \Ical} \to 0$. Hence, there exists $k_s \in \Ical$ such that for all $k \in \Ical$ with $k \geq k_s$ one finds $\|s_k\| \leq s_{thresh}$. Therefore, for all $k\in\Ical$ with $k \geq k_s$, it follows that $\rho_k \geq \eta$, as claimed.
\eproof

Next, we prove that the algorithm produces infinitely many accepted steps.

\blemma \label{lem.infinite_accepted}
  \textit{
  It holds that $|\Acal| = \infty$ and $\{s_k\}_{k\in\Acal}\to 0$.
  }
\elemma
\bproof  
  To derive a contradiction, suppose that $|\Acal| < \infty$.  This implies that there exists $k_0$ such that, for all $k \geq k_0$, one has $k \in \Rcal$ and $(x_k, g_k, H_k) = (x_{k_0}, g_{k_0}, H_{k_0})$.  From this fact and Assumption~\ref{ass.global}, it follows that $\|g_k\| \geq \epsilon$ for all $k \geq k_0$ for some $\epsilon \in \R{}_{++}$. From the fact that $k \in \Rcal$ for all $k \geq k_0$ and Lemma~\ref{lem.ratio_diverges}, it follows that $\{\lambda_k / \|s_k\|\} \to \infty$. This fact and $\|g_k\| \geq \epsilon$ for all $k \geq k_0$ imply that all the conditions of Lemma~\ref{lem.eventuallyaccepted} are satisfied for $\Ical := \{k \in \N{}_+ : k \geq k_0\}$; therefore, Lemma~\ref{lem.eventuallyaccepted} implies that for all sufficiently large $k \in \Ical$, one finds $\rho_k \geq \eta$ so that $k\in\Acal$, a contradiction.

To complete the proof, notice that the objective function values are monotonically decreasing.  Combining this with the condition in Step~\ref{step.accept}, the fact that $f$ is bounded below by $f_{\inf}$ (see Assumption~\ref{ass.f}), and $|\Acal| = \infty$, one deduces that $\{s_k\}_{k \in \Acal} \to 0$, as claimed.
\eproof

We now prove that there exists an infinite subsequence of iterates such that the sequence of gradients computed at those points converges to zero.

\blemma \label{lem.liminf}
  \textit{
  It holds that
  \bequationn
   \liminf_{k\in\N{}_+,k\to\infty} \|g_k\| = 0.
  \eequationn
  }
\elemma
\bproof
  To derive a contradiction, suppose that $\liminf_{k\in\N{}_+,k\to\infty} \|g_k\| > 0$, which along with the fact that $g_{k+1} = g_k$ for any $k \in \N{}_+\setminus\Acal$ means $\liminf_{k\in\Acal,k\to\infty} \|g_k\| > 0$. Thus, there exists $\epsilon \in \R{}_{++}$ such that
  \bequation\label{eq.g_lower}
    \|g_k\| \geq \epsilon\ \ \text{for all sufficiently large}\ \ k \in \Acal.
  \eequation
  
  Under \eqref{eq.g_lower}, let us prove that $\{\lambda_k\}_{k \in \Acal} \to \infty$. To derive a contradiction, suppose there exists an infinite $\Acal_\lambda \subseteq \Acal$ such that $\lambda_k \leq \lambda_{max}$ for some $\lambda_{max} \in \R{}_{++}$.  On the other hand, by $\{s_k\}_{k \in \Acal} \to 0$ (see Lemma~\ref{lem.infinite_accepted}) and \eqref{eq.residual_error}, it follows that $\{g_k + (H_k + \lambda_k I) s_k\}_{k \in \Acal_{\lambda}} \to 0$.  Combining the upper bound on $\{\lambda_k\}_{k\in\Acal_\lambda}$, the fact that $\{s_k\}_{k \in \Acal} \to 0$, and $\|H_k\| \leq H_{max}$ (see Lemma~\ref{lem.H_bounded}), it follows that $\{g_k\}_{k \in \Acal_\lambda} \to 0$, which violates \eqref{eq.g_lower}.  Therefore, $\{\lambda_k\}_{k\in\Acal} \to \infty$.

  Our next goal is to prove, still under \eqref{eq.g_lower}, that $k \in \Acal$ for all sufficiently large $k \in \N{}_+$.  To prove this, our strategy is to show that the sets of iterations involving a rejected step followed by an accepted step are finite.  In particular, let us define the index sets
  \begin{align*}
    \Rcal_1 &:= \{k\in\Rcal : \text{the condition in Step~\ref{step.reject.lambda} tests true and $(k+1)\in\Acal$} \} \ \ \text{and} \\
    \Rcal_2 &:= \{k\in\Rcal : \text{the condition in Step~\ref{step.reject.lambda} tests false and $(k+1)\in\Acal$} \}.
  \end{align*}
  We aim to prove that these are finite.  First, consider $\Rcal_1$.  To derive a contradiction, suppose that $|\Rcal_1| = \infty$.  By definition, for all $k\in\Rcal_1$, the condition in Step~\ref{step.reject.lambda} tests true, meaning $(s_{k+1},\lambda_{k+1})$ is found in Step~\ref{step.Pk} satisfying $\lambda_{k+1} / \|s_{k+1}\| \leq \sigmamax$.  On the other hand, since $(k+1) \in \Acal$ for all $k \in \Rcal_1$, it follows from Lemma~\ref{lem.infinite_accepted} that $\{s_{k+1}\}_{k \in \Rcal_1} \to 0$.  Combining the conclusions of these last two sentences shows that $\{\lambda_{k+1}\}_{k \in \Rcal_1} \to 0$.  However, this contradicts the conclusion of the previous paragraph, which showed that $\{\lambda_k\}_{k\in\Acal} \to \infty$.  Hence, we may conclude that $|\Rcal_1| < \infty$.  Now consider $\Rcal_2$.  To derive a contradiction, suppose that $|\Rcal_2| = \infty$.  The fact that the condition in Step~\ref{step.reject.lambda} tests false for $k\in\Rcal_2$ implies that $(s_{k+1},\lambda_{k+1})$ is found in Step~\ref{step.Pk} satisfying $\lambda_{k+1} / \|s_{k+1}\| \leq \gammau\lambda_k / \|s_k\|$. However, since $\{s_{k+1}\}_{k \in \Rcal_2} \to 0$ (see Lemma~\ref{lem.infinite_accepted}) and $\{\lambda_{k+1}\}_{k \in \Rcal_2} \to \infty$ (established in the previous paragraph), it follows that $\{\lambda_{k+1} / \|s_{k+1}\|\}_{k \in \Rcal_2} \to \infty$, which combined with the previously established inequality $\lambda_{k+1} / \|s_{k+1}\| \leq \gammau\lambda_k / \|s_k\|$ shows that $\{\lambda_k / \|s_k\|\}_{k \in \Rcal_2} \to \infty$.  Therefore, with~\eqref{eq.g_lower}, the conditions in Lemma~\ref{lem.eventuallyaccepted} hold for $\Iscr = \Rcal_2$, meaning that, for all sufficiently large $k \in \Rcal_2$, the inequality $\rho_k \geq \eta$ holds.  This contradicts the fact that $\Rcal_2 \subseteq \Rcal$; hence, we conclude that $\Rcal_2$ is finite.  Since $\Rcal_1$ and $\Rcal_2$ are finite, it follows from the logic of Algorithm~\ref{alg.inexact_regularized_Newton} that either $k\in\Acal$ for all sufficiently large $k$ or $k\in\Rcal$ for all sufficiently large~$k$.  By Lemma~\ref{lem.infinite_accepted}, it follows that $k\in\Acal$ for all sufficiently large $k$.

Thus far, we have proved under \eqref{eq.g_lower} that $\{\lambda_k\}_{k \in \Acal} \to \infty$ and that $k\in\Acal$ for all large $k \in \N{}_+$.  From this latter fact, it follows that there exists $k_\sigma$ such that $\sigmau_k = \sigmau_{k_\sigma} \in \R{}_{++}$ for all $k \geq k_\sigma$.  In addition, from Step~\ref{step.Pk}, it follows that for $k \geq k_\sigma$ one finds $\lambda_k / \|s_k\| \leq \sigmau_k = \sigmau_{k_\sigma} < \infty$.  However, this leads to a contradiction to the facts that $\{\lambda_k\}_{k \in \Acal} \to \infty$ and $\{s_k\}_{k \in \Acal} \to 0$ (see Lemma~\ref{lem.infinite_accepted}).  Overall, we have shown that \eqref{eq.g_lower} cannot be true, which proves the desired result.
\eproof

We close with our main global convergence result of this subsection, the proof of which borrows much from that of Theorem~3.14 in \cite{CurtRobiSama17}.

\btheorem\label{th.global}
  \textit{
  Under Assumptions~\ref{ass.subproblem_accuracy}, \ref{ass.f}, and \ref{ass.global}, it follows that
  \bequation\label{eq.lim}
    \lim_{k\in\N{}_+,k\to\infty} \|g_k\| = 0.
  \eequation
  }
\etheorem
\bproof
  For the purpose of reaching a contradiction, suppose that~\eqref{eq.lim} does not hold.  Combining this with the fact that $|\Acal| = \infty$ (see Lemma~\ref{lem.infinite_accepted}), it follows that there exists an infinite subsequence $\{t_i\} \subseteq \Acal$ (indexed over $i \in \N{}_+$) and a scalar $\epsilon > 0$ such that, for all $i \in \N{}_+$, one finds $\|g_{t_i}\| \geq 2 \epsilon > 0$.  Also, the fact that $|\Acal| = \infty$ and Lemma~\ref{lem.liminf} imply that there exists an infinite subsequence $\{\ell_i\} \subseteq \Acal$ (indexed over $i \in \N{}_+$) such that, for all $i \in \N{}_+$ and $k \in \N{}_+$ with $t_i \leq k < \ell_i$, one finds
  \bequation\label{eq.sandwich}
    \|g_k\| \geq \epsilon\ \ \text{and}\ \ \|g_{\ell_i}\| < \epsilon.
  \eequation
  Let us now restrict our attention to indices in the infinite index set
  \bequationn
    \Kcal := \{k \in \Acal : t_i \leq k < \ell_i\ \text{for some}\ i \in \N{}_+\}.
  \eequationn
  Observe from~\eqref{eq.sandwich} that, for all $k \in \Kcal$, it follows that $\|g_k\| \geq \epsilon$.  Also, from the definition of $\Acal$,
  \bequation\label{eq.f_decrease}
    f_k - f_{k+1} \geq \eta \|s_k\|^3\ \ \text{for all}\ \ k \in \Kcal \subseteq \Acal.
  \eequation
  Since $\{f_k\}$ is monotonically decreasing and bounded below, one finds that $\{f_k\} \to \underline{f}$ for some $\underline{f} \in \R{}$, which when combined with \eqref{eq.f_decrease} shows that
  \bequation\label{eq.delta_over_K}
    \lim_{k\in\Kcal,k\to\infty} \|s_k\| = 0.
  \eequation
  Using this fact, Lemma~\ref{lem.Cauchydecrease}, Assumption~\ref{ass.global}, and the Mean Value Theorem (as it is used in the proof of Lemma~\ref{lem.eventuallyaccepted} to yield \eqref{eq.mvt}), it follows that for all sufficiently large $k\in\Kcal$ one has
  \bequalin
    f_k - f_{k+1}
       =&\ f_k - q_k(s_k) + q_k(s_k) - f(x_k+s_k)     \\
    \geq&\ \frac{\|g_k\|}{6\sqrt{2}} \min \left\{\frac{\|g_k\|}{1+H_{max}},\frac{\|s_k\|}{\sqrt{6}}\sqrt{\frac{\|g_k\|}{2\|g_k\|+ \|s_k\|(\tfrac32 H_{max} + \kappa_1)}}\right\} - (g_{Lip} + \thalf H_{max}) \|s_k\|^2 \\
    \geq&\ \frac{\eps}{6\sqrt{2}} \frac{\|s_k\|}{\sqrt{6}}\sqrt{\frac{\|g_k\|}{2\|g_k\|+ \|s_k\|(\tfrac32 H_{max} + \kappa_1)}} 
       - (g_{Lip} + \thalf H_{max}) \|s_k\|^2.
  \eequalin
  It now follows from \eqref{eq.sandwich} and~\eqref{eq.delta_over_K} that, as $k \to \infty$ over $k\in\Kcal$, the square root term in the previous inequality converges to $1/\sqrt{2}$.  Since the second term in the previous inequality is of order $\|s_k\|^2$, the first term is of order $\|s_k\|$, and $1/\sqrt{2} > 1/\sqrt{3}$, one can thus conclude that $f_k - f_{k+1} \geq \eps\|s_k\|/36$ for all sufficiently large $k\in\Kcal$.  Consequently, it follows that for all sufficiently large $i \in \N{}_+$ one finds
  \bequalin
    \|x_{t_i} - x_{\ell_i}\|
      &\leq \sum_{k\in\Kcal,k=t_i}^{\ell_i-1} \|x_k - x_{k+1}\| \\
      &= \sum_{k\in\Kcal,k=t_i}^{\ell_i-1} \|s_k\| \leq \sum_{k\in\Kcal,k=t_i}^{\ell_i-1} \tfrac{36}{\epsilon}(f_k - f_{k+1}) = \tfrac{36}{\epsilon} (f_{t_i} - f_{\ell_i}).
  \eequalin
  Since $\{f_{t_i} - f_{\ell_i}\} \to 0$ (recall that $\{f_k\} \to \underline{f}$ monotonically) this implies that $\{\|x_{t_i} - x_{\ell_i}\|\} \to 0$, which, in turn, implies that $\{\|g_{t_i}-g_{\ell_i}\|\} \to 0$ because of the continuity of $g$.  However, this is a contradiction since, for any $i \in \N{}_+$, we have $\|g_{t_i} - g_{\ell_i}\| \geq \epsilon$ by the definitions of $\{t_i\}$ and $\{\ell_i\}$.  Overall, we conclude that our initial supposition must be false, implying that \eqref{eq.lim} holds.  
\eproof

\subsection{First-Order Complexity}\label{sec.complexity}

Our next goal is to prove, with respect to a prescribed positive threshold, a worst-case upper bound on the number of iterations required for our algorithm to reduce the norm of the gradient below the threshold.  In this subsection, along with Assumptions~\ref{ass.subproblem_accuracy}, \ref{ass.f}, and \ref{ass.global}, we add the following.

\bassumption\label{ass.complexity}
  \textit{
  The Hessian function $H$ is Lipschitz continuous on a path defined by the sequence of iterates and trial steps; in particular, it is Lipschitz continuous with a scalar Lipschitz constant $H_{Lip} > 0$ on the set $\{x_k + \tau s_k : k \in \N{}_+, \tau \in [0,1]\}$.
  }
\eassumption

We begin our analysis in this subsection by providing a lemma that shows that successful steps always result if $\lambda_k$ is sufficiently large relative to the size of the step.

\blemma \label{lem.lambda_large}
  \textit{
  For any $k \in \N{}_+$, if the pair $(s_k,\lambda_k)$ satisfies
  \bequation\label{eq.lambda_large}
    \lambda_k \geq (H_{Lip} + \kappa_2 + 2\eta)\|s_k\|,
  \eequation
  then $\rho_k \geq \eta$.
  }
\elemma
\bproof
  It follows from Assumption~\ref{ass.complexity} and Taylor's expansion with Lagrange remainder that there exists $\xbar_k$ on the line segment $[x_k,x_k+s_k]$ such that 
  \bequation \label{eq.modeldeviation}
    q_k(s_k) - f(x_k + s_k) = \thalf s_k^T (H_k - H(\xbar_k))s_k \geq - \thalf H_{Lip} \|s_k\|^3.
  \eequation
  Also, it follows from \eqref{eq.subspace_optimality} that 
  \bequation \label{eq.modeldecrease}
    \baligned
      f_k - q_k(s_k) 
      &= - g_k^T s_k - \thalf s_k^T H_k s_k \\
      &= - s_k^T (g_k + (H_k + \lambda_k I)s_k) + \thalf \lambda_k \|s_k\|^2 + \thalf s_k^T (H_k + \lambda_k I)s_k \\
      &\geq - \thalf s_k^T (H_k + \lambda_k I)s_k - \thalf \kappa_2 \|s_k\|^3 + \thalf \lambda_k \|s_k\|^2 + \thalf s_k^T (H_k + \lambda_k I)s_k \\
      &= - \thalf \kappa_2 \|s_k\|^3 + \thalf \lambda_k \|s_k\|^2.
    \ealigned
  \eequation
  From \eqref{eq.modeldeviation} and \eqref{eq.modeldecrease}, it follows that
  \bequationn
    \baligned
      f_k - f(x_k + s_k) 
      &= f_k - q_k(s_k) + q_k(s_k) - f(x_k+s_k) \\
      &\geq \thalf \lambda_k \|s_k\|^2 - \thalf \kappa_2 \|s_k\|^3 - \thalf H_{Lip} \|s_k\|^3,
    \ealigned 
  \eequationn
which together with~\eqref{eq.lambda_large} implies that $\rho_k \geq \eta$, as claimed.
\eproof

We now prove that the sequence $\{\sigmau_k\}$ is bounded above.

\blemma\label{lem.sigma_bounded}
  \textit{
  There exists a scalar constant $\sigma_{max} \in \R{}_{++}$ such that, for all $k \in \N{}_+$,
  \bequationn
    \sigmau_k \leq \sigma_{max}.
  \eequationn
  }
\elemma
\bproof
  Consider any $k \in \N{}_+$.  If $s_k$ is accepted (i.e., $k\in\Acal$), then $\sigmau_{k+1} \gets \sigmau_k$.  On the other hand, if $s_k$ is rejected (i.e., $k\in\Rcal$), then it follows from Step~\ref{step.reject.lambda_small} and Step~\ref{step.reject.lambda_large} that $\sigmau_{k+1} \leq \max\{\sigmamax,\gammau \lambda_k/\|s_k\|\}$.  Moreover, since $k\in\Rcal$, meaning that $\rho_k < \eta$, it follows from Lemma~\ref{lem.lambda_large} that $\lambda_k / \|s_k\|$ is bounded above by $(H_{Lip} + \kappa_2 + 2 \eta)$. Thus, it follows that $\sigmau_{k+1} \leq  \max\{\sigmamax,\gammau(H_{Lip} + \kappa_2 + 2 \eta) \}$ for all $k \in\Rcal$.  Overall, the desired result follows for any $\sigma_{max} \geq \max\{\sigmamax,\gammau(H_{Lip} + \kappa_2 + 2 \eta)\}$.
\eproof

We now establish a lower bound on the norm of any accepted trial step.  

\blemma \label{lem.s-and-g}
  \textit{
  For all $k \in \Acal$, it follows that
  \bequationn
    \|s_k\| \geq \(\thalf H_{Lip} + 2\sigma_{max} + \kappa_3\)^{-1/2} \|g_{k+1}\|^{1/2}.
  \eequationn
  }
\elemma
\bproof
  Let $k \in \Acal$. It follows that
  \begin{align}
    \|g_{k+1}\|
      &\leq \| g_{k+1} - (g_k + (H_k + \lambda_k I)s_k) \| + \| g_k + (H_k + \lambda_k I)s_k \| \nonumber \\
      &\leq \| g_{k+1} - (g_k + H_ks_k) \| + \lambda_k \|s_k\| + \| g_k + (H_k + \lambda_k I)s_k \|. \label{eq.mvt_new}
  \end{align}
  By Taylor's theorem and Assumption~\ref{ass.complexity}, the first term on the right-hand side of this inequality satisfies
  \bequationn
    \baligned
      \| g_{k+1} - (g_k + H_ks_k) \|
        &\leq \left\| \int_0^1 (H(x_k + \tau s_k) - H_k)s_k d\tau \right\| \\
        &\leq \int_0^1 \|H(x_k + \tau s_k) - H_k\| d\tau \cdot \|s_k\| \\
        &\leq \int_0^1 \tau d\tau \cdot H_{Lip} \|s_k\|^2 = \thalf H_{Lip} \|s_k\|^2.
    \ealigned
  \eequationn
  Combining this with \eqref{eq.mvt_new} and observing Step~\ref{step.Pk}, \eqref{eq.residual_error}, and Lemma~\ref{lem.sigma_bounded}, it follows that
  \bequalin
    \|g_{k+1}\|
      &\leq \thalf H_{Lip} \|s_k\|^2 + 2\frac{\lambda_k}{\|s_k\|} \|s_k\|^2 + \kappa_3 \|s_k\|^2 \\
      &\leq \thalf H_{Lip} \|s_k\|^2 + 2\sigma_{max} \|s_k\|^2 + \kappa_3 \|s_k\|^2,
  \eequalin
  which, after rearrangement, completes the proof.  
\eproof

We are now prepared to prove a worst-case upper bound on the total number of accepted steps that may occur for iterations in which the norm of the gradient of the objective is above a positive threshold.

\blemma\label{lem.num.successful}
  \textit{
  For any $\epsilon\in\R{}_{++}$, the total number of elements in the index set
  \bequationn
    \Kcal_\epsilon := \{k \in \N{}_+ : k \geq 1,\ (k-1) \in \Acal,\ \|g_k\| > \epsilon\}
  \eequationn
  is at most
  \bequation\label{eq.iter_limit}
    \left\lfloor \(\frac{f_0 - f_{\inf}}{\eta (\thalf H_{Lip} + 2\sigma_{max} + \kappa_3)^{-3/2}}\) \epsilon^{-3/2} \right\rfloor =: N_\Acal(\epsilon) \geq 0.
  \eequation
  }
\elemma
\bproof
  The proof follows in a similar manner as that of Lemma 3.20 in \cite{CurtRobiSama17}.  By Lemma~\ref{lem.s-and-g}, it follows that, for all $k \in \Kcal_\epsilon$, one finds
  \bequalin
    f_{k-1} - f_k
      &\geq \eta \|s_{k-1}\|^3 \\
      &\geq \eta (\thalf H_{Lip} + 2\sigma_{max} + \kappa_3)^{-3/2} \|g_k\|^{3/2} \\
      &\geq \eta (\thalf H_{Lip} + 2\sigma_{max} + \kappa_3)^{-3/2} \epsilon^{3/2}.
  \eequalin
  In addition, it follows from Theorem~\ref{th.global} that $|\Kcal_\epsilon| < \infty$.  Hence, the reduction in $f$ obtained up to the largest index in $\Kcal_\epsilon$, call it $\kbar_\epsilon$, satisfies
  \bequationn
    f_0 - f_{\kbar_\epsilon} = \sum_{k=1}^{\kbar_\epsilon} (f_{k-1} - f_k) \geq \sum_{k\in\Kcal_\epsilon} (f_{k-1} - f_k) \geq |\Kcal_\epsilon| \eta (\thalf H_{Lip} + 2\sigma_{max} + \kappa_3)^{-3/2} \epsilon^{3/2}.
  \eequationn
  Rearranging this inequality to yield an upper bound for $|\Kcal_\epsilon|$ and using the fact that $f_0 - f_{\inf} \geq f_0 - f_{\kbar_\epsilon}$, one obtains the desired result.  
\eproof

In order to prove a result similar to Lemma~\ref{lem.num.successful} for the \emph{total} number of iterations with $\|g_k\| > \epsilon$, we require an upper bound on the total number of trial steps that may be rejected between accepted steps. To this end, let us define, for a given $\khat \in \Acal \cup \{0\}$, the iteration number and corresponding set
\bequalin
  k_\Acal(\khat) &:= \min\{k \in \Acal : k > \khat\} \\ \text{and}\ \ 
  \Ical(\khat)   &:= \{k \in \N{}_+: \khat < k < k_\Acal(\khat)\},
\eequalin
i.e, we let $k_\Acal(\khat)$ be the smallest of all iteration numbers in $\Acal$ that is strictly larger than~$\khat$, and we let $\Ical(\khat)$ be the set of iteration numbers between $\khat$ and $k_\Acal(\khat)$.

We now show that the number of rejected steps between the first iteration and the first accepted step, or between consecutive accepted steps, is bounded above.

\blemma\label{lem.Cbounded}
  \textit{
  For any $\khat \in \Acal \cup \{0\}$, it follows that
  \bequationn
    |\Ical(\khat)| \leq 1 + \left \lfloor \frac{1}{ \log(\gammal)} \log\left(\frac{\sigma_{max}}{\sigmamin}\right) \right \rfloor =: N_\Rcal \geq 0.
  \eequationn
  }
\elemma
\bproof
  The proof follows in a similar manner as for Lemma 3.24 in \cite{CurtRobiSama17}.  First, the result holds trivially if $|\Ical(\khat)| = 0$.  Thus, we may assume that $|\Ical(\khat)| \geq 1$.  Since $(\khat +1) \in \Rcal$ by construction, it follows from Steps~\ref{step.reject.lambda}--\ref{step.reject.lambda_large} and Step~\ref{step.Pk} that $\lambda_{\khat+2}/\|s_{\khat+2}\| \geq \sigmamin$, which, due to the lower bound on $\lambda_{k+1}/\|s_{k+1}\|$ in Step~\ref{step.reject.lambda_large} and Step~\ref{step.Pk}, leads to
  \bequationn
    \lambda_{k_\Acal(\khat)} \geq \sigmamin\(\gammal\)^{k_\Acal(\khat)-\khat-2} \|s_{k_\Acal(\khat)}\|.
  \eequationn
  Combining this with Step~\ref{step.Pk} and Lemma~\ref{lem.sigma_bounded}
  shows that
  \bequationn
    \sigma_{max} 
    \geq \sigmau_{k_\Acal(\khat)}
    \geq \lambda_{k_\Acal(\khat)}/\|s_{k_\Acal(\khat)}\| 
    \geq \sigmamin\(\gammal\)^{k_\Acal(\khat)-\khat-2}.
  \eequationn
  After rearrangement, it now follows that
  \bequationn
    k_\Acal(\khat) - \khat -2 \leq \frac{1}{\log(\gammal)}\log\left( \frac{\sigma_{max}}{\sigmamin} \right).
  \eequationn
  The desired result follows from this inequality since $|\Ical(\khat)| = k_\Acal(\khat) - \khat - 1$.  
\eproof

We are now prepared to prove our main complexity result of this subsection.

\btheorem\label{th.main_complexity}
  \textit{
  Under Assumptions~\ref{ass.subproblem_accuracy}, \ref{ass.f}, \ref{ass.global}, and \ref{ass.complexity}, for a scalar $\epsilon\in\R{}_{++}$, the total number of elements in the index set $\{k \in \N{}_+ : \|g_k\| > \epsilon\}$ is at most
  \bequation\label{eq.first_order_bound}
    N(\epsilon) := 1 + N_\Rcal N_\Acal(\epsilon),
  \eequation
  where $N_\Acal(\epsilon)$ and $N_\Rcal$ are defined in Lemmas~\ref{lem.num.successful} and \ref{lem.Cbounded}, respectively.  Consequently, for any $\overline\epsilon \in \R{}_{++}$, it follows that $N(\epsilon) = \Ocal(\epsilon^{-3/2})$ for all $\epsilon \in (0,\overline\epsilon]$.
  }
\etheorem
\bproof
  Without loss of generality, we may assume that at least one iteration is performed.  Lemma~\ref{lem.num.successful} guarantees that the total number of elements in the index set $\{k \in \Acal : k \geq 1,\ \|g_k\| > \epsilon\}$ is at most $N_\Acal(\epsilon)$, where, immediately prior to each of the corresponding accepted steps, Lemma~\ref{lem.Cbounded} guarantees that at most $N_\Rcal$ trial steps are rejected.  Accounting for the first iteration, the desired result follows.  
\eproof

\subsection{Second-Order Global Convergence and Complexity}\label{sec.second_order}

Our goal in this subsection is to prove results showing that, in some sense, the algorithm converges to second-order stationarity and does so with a worst-case iteration complexity on par with the methods in \cite{CartGoulToin11b} and \cite{CurtRobiSama17}.  In particular, our results show that if the algorithm computes each search direction to satisfy a curvature condition over a subspace, then second-order stationarity is reached in a manner that depends on the subspaces.

In this subsection, we make the following additional assumption about the subproblem solver.

\bassumption\label{ass.2nd}
  \textit{
  For all $k\in\N{}_+$, let $\Lcal_k \subseteq \R{n}$ denote a subspace with an orthonormal basis formed from the columns of a matrix $R_k$.  The step $s_k$ satisfies
  \bequation\label{eq.2nd}
    \xi(R_k^T H_k R_k) \geq - \kappa_4\|s_k\|
  \eequation
  for some $\kappa_4\in \R{}_+$, where $\xi(R_k^T H_k R_k)$ indicates the smallest eigenvalue of $R_k^T H_k R_k$.
  }
\eassumption

\noindent
This assumption is reasonable, e.g., in cases when $s_k$ is computed by solving problem~\ref{prob.primal-dual} with the component $s$ restricted to a subspace of $\R{n}$.  We refer the reader to Theorem~\ref{th.consistency-subspace} for a proof of this fact, which also reveals that this assumption is congruous with Assumption~\ref{ass.subproblem_accuracy}.

Under this assumption, we have the following second-order convergence result.

\btheorem\label{th.2nd_global}
  \textit{
  Suppose Assumptions~\ref{ass.subproblem_accuracy}, \ref{ass.f}, \ref{ass.global}, \ref{ass.complexity}, and~\ref{ass.2nd} hold.  It follows that
  \bequationn
    \liminf_{k \in \Acal, k \to \infty} \xi(R_k^T H_k R_k) \geq 0.
  \eequationn
  }
\etheorem
\bproof
  The result follows from \eqref{eq.2nd} since $\{s_k\}_{k \in \Acal} \to 0$ (see Lemma~\ref{lem.infinite_accepted}).
\eproof

As a consequence of Theorem~\ref{th.2nd_global}, if the sequence $\{R_k\}_{k\in\Acal}$ tends toward full-dimensionality as $k \to \infty$, then any limit point $x_*$ of $\{x_k\}$ must have $H(x_*) \succeq 0$.

Our next goal is to prove a worst-case iteration complexity result for achieving second-order stationarity in a sense similar to that in Theorem~\ref{th.2nd_global}.  Toward this end, we first prove the following lemma, which is similar to Lemma~\ref{lem.num.successful}.

\blemma\label{lem.num.successful.2nd}
  \textit{
  For any $\epsilon\in\R{}_{++}$, the total number of elements in the index set
  \bequationn
    \Kcal_{\epsilon,\xi} := \{k \in \N{}_+ : k \geq 1,\ (k-1) \in \Acal,\ \xi(R_k^T H_k R_k)  < -\epsilon\}
  \eequationn
  is at most
  \bequation\label{eq.iter_limit.2nd}
    \left\lfloor \(\frac{f_0 - f_{\inf}}{\eta \kappa_4^{-3}}\) \epsilon^{-3} \right\rfloor =: N_{\Acal,\xi}(\epsilon) \geq 0.
  \eequation
  }
\elemma
\bproof
  Under Assumption~\ref{ass.2nd}, it follows that, for all $k \in \Kcal_{\epsilon,\xi}$, one finds
  \bequationn
    f_{k-1} - f_k \geq \eta \|s_{k-1}\|^3 \geq \eta \left(\frac{-\xi(R_k^T H_k R_k)}{\kappa_4}\right)^3 \geq \eta \kappa_4^{-3} \epsilon^3.
  \eequationn
  It follows from this inequality, the fact that $f$ is monotonically decreasing over the sequence of iterates, and Assumption~\ref{ass.f} that
  \bequationn
    f_0 - f_{\inf} \geq \sum_{k\in\Kcal_{\epsilon,\xi}} (f_{k-1} - f_k) \geq |\Kcal_{\epsilon,\xi}| \eta \kappa_4^{-3} \epsilon^3.
  \eequationn
  Rearranging this inequality to yield an upper bound for $|\Kcal_{\epsilon,\xi}|$ gives the result.  
\eproof

We close with the following second-order complexity result.

\btheorem\label{th.main_complexity.2nd}
  \textit{
  Under Assumptions~\ref{ass.subproblem_accuracy}, \ref{ass.f}, \ref{ass.global}, \ref{ass.complexity}, and \ref{ass.2nd},  for any pair of scalars $(\epsilon_1,\epsilon_2) \in \R{}_{++} \times \R{}_{++}$, the number of elements in the index set
  \bequationn
    \{k \in \N{}_+ : \|g_k\| > \epsilon_1\ \vee\ \xi(R_k^T H_k R_k) < -\epsilon_2\}
  \eequationn
  is at most
  \bequation\label{eq.second_order_bound}
    N(\epsilon_1,\epsilon_2) := 1 + N_\Rcal \max\{N_\Acal(\epsilon_1),N_{\Acal,\xi}(\epsilon_2)\},
  \eequation
  where $N_\Acal(\cdot)$, $N_\Rcal$, and $N_{\Acal,\xi}(\cdot)$ are defined in Lemmas~\ref{lem.num.successful}, \ref{lem.Cbounded}, and \ref{lem.num.successful.2nd}, respectively.  Consequently, for any pair of scalars $(\overline\epsilon_1,\overline\epsilon_2) \in \R{}_{++} \times \R{}_{++}$, it follows that
  \bequationn
    N(\epsilon_1,\epsilon_2) = \Ocal(\max\{\epsilon_1^{-3/2},\epsilon_2^{-3}\})\ \ \text{for all}\ \ (\epsilon_1,\epsilon_2) \in (0,\overline\epsilon_1] \times (0,\overline\epsilon_2].
  \eequationn
  }
\etheorem
\bproof
  The proof follows in a similar manner as that of Theorem~\ref{th.main_complexity} by additionally incorporating the bound proved in Lemma~\ref{lem.num.successful.2nd}.
\eproof

\section{Algorithm Instances}\label{sec.instances}

Algorithm~\ref{alg.inexact_regularized_Newton} is a broad framework containing, amongst other algorithms, \ARC{} and \TRACE.  Indeed, the proposed framework and its supporting analyses cover a wide range of algorithms as long as the pairs in the sequence $\{(s_k,\lambda_k)\}$ satisfy Assumption~\ref{ass.subproblem_accuracy}.

In this section, we show that \ARC{} and \TRACE{} are special cases of our proposed framework in that the steps these algorithms accept would also be acceptable for our framework, and that the procedures followed by these methods after a step is rejected are consistent with our framework.  We then introduce an instance of our frameowork that is new to the literature.  (If desired for the guarantees in \S\ref{sec.second_order}, one could also mind whether the elements in the sequence $\{(s_k,\lambda_k)\}$ satisfy Assumption~\ref{ass.2nd}.  However, for brevity in this section, let us suppose that one is interested only in Assumption~\ref{ass.subproblem_accuracy}.)

\subsection{\ARC{} as a Special Case}\label{subsec.arc}

The \ARC{} method, which was inspired by the work in \cite{Grie81} and \cite{NestPoly06}, was first proposed and analyzed in \cite{CartGoulToin11a,CartGoulToin11b}.  In these papers, various sets of step computation conditions are considered involving exact and inexact subproblem solutions yielding different types of convergence and worst-case complexity guarantees.  For our purposes here, we consider the more recent variant of \ARC{} stated and analyzed as ``\textsc{ar}$p$'' with $p=2$ in \cite{BirgGardMartSantToin17}.  (For ease of comparison, we consider this algorithm when their regularization parameter update---see Step~4 in their algorithm---uses $\eta_1 = \eta_2$.  Our algorithm is easily extended to employ a two-tier acceptance condition, involving two thresholds $\eta_1$ and $\eta_2$, as is used in \cite{BirgGardMartSantToin17} and \cite{CartGoulToin11a,CartGoulToin11b}.)

Suppose that a trial step $s_k$ is computed by this version of \ARC{}.  In particular, let us make the reasonable assumption that the subproblem for which~$s_k$ is an approximate solution is defined by some regularization value $\sigma_k \in [\sigmal_k,\sigmau_k]$ (with $\sigmal_k \geq \sigma_{min}$ since \ARC{} ensures that $\sigma_k \geq \sigma_{min} \in \R{}_{++}$ for all $k \in \N{}$) and that this subproblem is minimized over a subspace $\Lcal_k$ such that $g_k \in \Lcal_k$ (see Appendix~\ref{app.subproblem.reduced}).  As is shown using a similar argument as in the proof of our Theorem~\ref{th.consistency-subspace}$(b)$, one can show under these conditions that $(s_k,\lambda_k)$ with $\lambda_k = \sigma_k\|s_k\|$ satisfies \eqref{eq.Cauchy_decrease}.  In addition, considering the algorithm statement in \cite{BirgGardMartSantToin17}, but using our notation, one is required to have
\bequationn
  g_k^Ts_k + \thalf s_k^TH_ks_k + \lambda_k \|s_k\|^2 < 0\ \ \text{and}\ \ \|g_k + (H_k + \lambda_kI)s_k\| \leq \theta \|s_k\|^2\ \ \text{for some}\ \ \theta \in \R{}_{++}.
\eequationn
It is easily seen that $(s_k,\lambda_k)$ satisfying these conditions also satisfies \eqref{eq.subspace_optimality}--\eqref{eq.residual_error} for any $(\kappa_1,\kappa_2,\kappa_3)$ such that $\kappa_1 \geq \thalf H_{max}$ and $\kappa_3 \geq \theta$.  Overall, we have shown that a trial step $s_k$ computed by this version of \ARC{} satisfies Assumption~\ref{ass.subproblem_accuracy}, meaning that it satisfies the condition in Step~\ref{step.Pk} in Algorithm~\ref{alg.inexact_regularized_Newton}.  If this trial step is accepted by \ARC{}, then this means that $f_k - f(x_k + s_k) \geq \eta_1(f_k - q_k(s_k))$.  Along with \cite[Lemma~2.1]{BirgGardMartSantToin17}, this implies that $f_k - f(x_k + s_k) \geq \tfrac13 \eta \sigma_k \|s_k\|^3$, meaning that $\rho_k \geq \tfrac13 \eta_1 \sigma_{min}$.  Hence, this trial step would also be accepted in Algorithm~\ref{alg.inexact_regularized_Newton} under the assumption that $\eta \in (0,\tfrac13 \eta_1 \sigma_{min}]$.

Finally, if a trial step is rejected in this version of \ARC{}, then $\sigma_{k+1}$ is set to a positive multiple of $\sigma_k$.  This is consistent with the procedure after a step rejection in Algorithm~\ref{alg.inexact_regularized_Newton}, where it is clear that, with appropriate parameter choices, one would find $\sigma_{k+1} \in [\sigmal_{k+1},\sigmau_{k+1}]$.

\subsection{\TRACE{} as a Special Case}\label{subsec.trace}


\TRACE{} is proposed and analyzed in \cite{CurtRobiSama17}.  Our goal in this subsection is to show that, with certain parameter settings, a trial step that is computed and accepted by \TRACE{} could also be one that is computed and accepted by Algorithm~\ref{alg.inexact_regularized_Newton}, and that the procedures for rejecting a step in \TRACE{} are consistent with those in Algorithm~\ref{alg.inexact_regularized_Newton}.  Amongst other procedures, \TRACE{} involves dynamic updates for two sequences, $\{\delta_k\}$ and $\{\Delta_k\}$.  The elements of $\{\delta_k\}$ are the trust region radii while $\{\Delta_k\}$ is a monotonically nondecreasing sequence of upper bounds for the trust region radii; consequently, $\|s_k\| \leq \delta_k \leq \Delta_k$ with $\Delta_{k+1} \geq \Delta_k$ for all $k \in \N{}$.  For simplicity in our discussion here, let us assume that $\|s_k\| < \Delta_k$ for all $k \in \N{}$.  This is a fair assumption since, as shown in \cite[Lemma~3.11]{CurtRobiSama17}, the manner in which $\{\Delta_k\}$ is set ensures that $\|s_k\| = \Delta_k$ only a finite number of times in any run.

In \TRACE{}, during iteration $k \in \N{}$, a trust region radius $\delta_k \in \R{}_{++}$ is given and a trial step $s_k$ and regularization value $\lambda_k$ are computed satisfying the standard trust region subproblem optimality conditions
\bequationn
  g_k + (H_k + \lambda_k I)s_k = 0,\ \ H_k + \lambda_k I \succeq 0,\ \ \text{and}\ \ \lambda_k (\delta_k - \|s_k\|) =0,\ \ \text{where}\ \ (\lambda_k,\delta_k-\|s_k\|) \geq 0.
\eequationn
By the first of these conditions, the pair $(s_k,\lambda_k)$ clearly satisfies \eqref{eq.subspace_optimality}--\eqref{eq.residual_error}.  In addition, one can use standard trust region theory, in particular related to Cauchy decrease (see \cite{ConnGoulToin00} or \cite{NoceWrig06}), to show that the pair also satisfies \eqref{eq.Cauchy_decrease}.  Overall, assuming that the pair $(\sigmal_k,\sigmau_k)$ is set such that $\lambda_k/\|s_k\| \in [\sigmal_k,\sigmau_k]$, it follows that Assumption~\ref{ass.subproblem_accuracy} is satisfied, meaning that \TRACE{} offers the condition in Step~\ref{step.Pk} in Algorithm~\ref{alg.inexact_regularized_Newton}.  If the trial step $s_k$ is subsequently accepted by \TRACE{}, then it would also be accepted by Algorithm~\ref{alg.inexact_regularized_Newton} since both algorithms use the same step acceptance condition.

Now suppose that a trial step is not accepted in \TRACE{}.  This can occur in two circumstances.  It can occur if $\rho_k \geq \eta$ while $\lambda_k > \sigma_k\|s_k\|$, in which case the trust region radius is \emph{expanded} and a new subproblem is solved.  By the proof of \cite[Lemma~3.7]{CurtRobiSama17}, the solution of this new subproblem yields (in iteration $k+1$ in \TRACE{}) the relationship that $\lambda_{k+1}/\|s_{k+1}\| \leq \sigma_{k+1} = \sigma_k$.  Hence, under the same assumption as above that the pair $(\sigmal_k,\sigmau_k)$ is set such that $\lambda_k/\|s_k\| \in [\sigmal_k,\sigmau_k]$, this shows that the procedure in \TRACE{} involving \emph{an expansion of the trust region radius and the computation of the subsequent trial step} yields a trial step that would be offered in a \emph{single iteration} in Algorithm~\ref{alg.inexact_regularized_Newton}.  The other circumstance in which a trial step is rejected in \TRACE{} is when $\rho_k < \eta$, in which case the trust region radius is contracted.  In this case, one can see that the outcome of the \textsc{contract} subroutine in \TRACE{} is consistent with Steps~\ref{step.reject.lambda}--\ref{step.reject.lambda_large} of Algorithm~\ref{alg.inexact_regularized_Newton} in the sense that the solution of the subsequent subproblem in \TRACE{} will have $\lambda_{k+1}/\|s_{k+1}\| \in [\sigmamin,\sigmamax]$ (if $\lambda_k < \sigmamin\|s_k\|$) or $\lambda_{k+1}/\|s_{k+1}\|$ within a range defined by positive multiples of $\lambda_k/\|s_k\|$; see Lemmas~3.17 and 3.23 in \cite{CurtRobiSama17}.


\subsection{A Hybrid Algorithm}\label{subsec.hybrid}

The primary distinguishing feature of our algorithm instance is the manner in which we compute the pair $(s_k,\lambda_k)$ in Step~\ref{step.Pk} of Algorithm~\ref{alg.inexact_regularized_Newton}.  Our newly proposed hybrid algorithm considers two cases.

\bitemize
  \item[\textbf{Case 1:}] $\sigmal_k > 0$. 
In this case, we find a pair $(s_k,\lambda_k)$ by solving problem~\eqref{prob.primal-dual.sigma-reduced} over a sequence of increasingly higher dimensional Krylov subspaces as described in~\cite{CartGoulToin11a} until~\eqref{eq.subproblem_accuracy} and~\eqref{eq.2nd} are satisfied.  The reason we know that~\eqref{eq.subproblem_accuracy} and~\eqref{eq.2nd} will eventually be satisfied can be seen as follows. Solving problem~\eqref{prob.primal-dual.sigma-reduced} over a Krylov subspace is equivalent to solving problem~\eqref{prob.primal-dual_subspace.sigma-reduced} with an appropriate choice of $R_k$ as a basis for that Krylov subspace, then setting $s_k = R_k v_k$. Then, it follows from Theorem~\ref{th.equivalency_subspace}$(i)$ that solving \eqref{prob.primal-dual_subspace.sigma-reduced} is equivalent to solving~\eqref{prob.primal-dual_subspace-reduced}, which in turn is equivalent to solving~\eqref{prob.primal-dual_subspace} in the sense that if $(v_k,\lambda_k,\betal_k,\betau_k,\betan_k)$ is a first-order primal-dual solution of problem~\eqref{prob.primal-dual_subspace-reduced}, then $(s_k,\lambda_k,\betal_k,\betau_k,\betan_k)$ with $s_k = R_kv_k$ is a  solution of problem~\eqref{prob.primal-dual_subspace}. Finally, we need only note from Theorem~\ref{th.consistency-subspace} that solutions to problem~\eqref{prob.primal-dual_subspace} satisfy~\eqref{eq.Cauchy_decrease} for all Krylov subspaces $\Lcal_k$ (recall that $g_k$ is contained in all Krylov subspaces), \eqref{eq.subspace_optimality} for all Krylov subspaces, \eqref{eq.residual_error} if the Krylov subspace $\Lcal_k$ includes enough of the space (in the worst case, $\Lcal_k = \R{n}$), and \eqref{eq.2nd} for all Krylov subspaces.
  \item[\textbf{Case 2:}] $\sigmal_k = 0$. In this case, we begin by applying the linear CG method in an attempt to solve the linear system $H_k s = -g_k$, which iteratively solves
\bequation\label{prob.quadratic}
  \min_{s\in\R{n}}\ q_k(s)
\eequation
over a sequence of expanding Krylov subspaces.   
One of two outcomes is possible.  First, the CG algorithm may ultimately identify a vector $s_k$ such that $(s_k,\lambda_k)$ with $\lambda_k = 0$ satisfies~\eqref{eq.subproblem_accuracy} and~\eqref{eq.2nd}.  Second, the CG algorithm may never identify a vector $s_k$ such that $(s_k,\lambda_k)$ with $\lambda_k = 0$ satisfies~\eqref{eq.subproblem_accuracy} and~\eqref{eq.2nd}.  Indeed, this might occur if CG encounters a direction of negative curvature---in which case we terminate CG immediately---or if CG solves \eqref{prob.quadratic} accurately or reaches an iteration limit, and yet at least one condition in \eqref{eq.subproblem_accuracy}/\eqref{eq.2nd} is not satisfied.  In such a case, we choose to reset $\sigmal_k \in (0,\sigmau_k]$, then solve problem~\eqref{prob.primal-dual.sigma-reduced} over a sequence of expanding Krylov subspaces as described in Case~1.  In this manner, we are guaranteed to identify a pair $(s_k,\lambda_k)$ satisfying~\eqref{eq.subproblem_accuracy} and~\eqref{eq.2nd} as required.
\eitemize

\section{Implementation and Numerical Results}\label{sec.numerical}

We implemented two algorithms in \texttt{MATLAB}, one following the strategy in \S\ref{subsec.hybrid} and, for comparison purposes, one following the \ARC{} algorithm in \cite{CartGoulToin11b} with ideas from \cite{BirgGardMartSantToin17}.  We refer to our implementation of the former as \iRNewton{}, for inexact Regularized Newton, and to our implementation of the latter as \iARC{}, for inexact \ARC{}.  In this section, we describe our approach for computing the pairs $\{(s_k,\lambda_k)\}$ in \iRNewton{} and \iARC{}, as well as other implementation details, and discuss the results of numerical experiments on a standard set of nonlinear optimization test problems.

\subsection{Implementation Details}\label{subsec.implementation}

Let us begin by noting that the implemented algorithms terminate in iteration $k \in \N{}_+$ if
\bequationn
  \|g_k\|_{\infty} \leq 10^{-6} \max\{\|g_0\|_{\infty},1\}.
\eequationn
We chose not to employ a termination test based on a second-order stationarity condition.  Correspondingly, neither of the algorithms check a second-order condition when computing a trial step; e.g., in \iRNewton{}, we are satisfied with a step satisfying \eqref{eq.subproblem_accuracy} and do not check \eqref{eq.2nd}.  In addition, for practical purposes, we set an maximum iteration limit of $10^6$, a time limit of four hours, and a minimum step norm limit of $10^{-20}$.  For reference, the input parameter values we used are given in Table~\ref{tbl.parameters}.  We chose these values as ones that worked well on our test set for both implemented algorithms.

\btable[ht]
  \setlength{\tabcolsep}{12pt}
  \centering
  \caption{Input parameters for \iARC{} and \iRNewton}
  \smallskip
  
  \label{tbl.parameters}
  \texttt{
  \btabular{|l|l||l|l||l|l||l|l|}
    \hline
    $\eta_1$ & 1.0e-16 & $\gamma_0$ & 2.0e-01 & $\kappa_1$ & 1.0e+00 & $\sigmamin$ & 1.0e-10 \\
    $\eta_2$ & 1.0e-01 & $\gammal$  & 1.0e+01 & $\kappa_2$ & 1.0e+00 & $\sigmamax$ & 1.0e+20 \\
             &         & $\gammau$  & 2.0e+02 & $\kappa_3$ & 1.0e+00 &             &         \\
    \hline
  \etabular
  }
\etable

For both implemented algorithms, we employ a sequence $\{\sigma_k\}$ that is updated dynamically.  In \iARC{}, this sequence is handled as described in \cite{CartGoulToin11b}, namely,
\bequationn
  \sigma_{k+1} \gets \bcases \max\{\sigmamin,\gamma_0 \sigma_k\} & \text{if $\frac{f_k - f(x_k + s_k)}{f_k - c_k(s_k;\sigma_k)} \geq \eta_2$} \\[5pt] \sigma_k & \text{if $\frac{f_k - f(x_k + s_k)}{f_k - c_k(s_k;\sigma_k)} \in [\eta_1,\eta_2)$} \\[5pt] \gammal \sigma_k & \text{if $\frac{f_k - f(x_k + s_k)}{f_k - c_k(s_k;\sigma_k)} < \eta_1$} \ecases
\eequationn
The value $\sigma_k$ is used in defining $c_k(\cdot;\sigma_k)$ (recall \eqref{eq.c}) that is minimized approximately to compute the trial step $s_k$ for all $k \in \N{}_+$.  In particular, the implementation iteratively constructs Krylov subspaces of increasing dimension using the Lanczos process, where for each subspace we employ the \texttt{RQS} function from the \texttt{GALAHAD} software library (see \cite{GoulOrbaToin03} and \cite{GoulRobiThor10}) to minimize $c_k(\cdot;\sigma_k)$ over the subspace.  If the subspace is full-dimensional or the resulting step $s_k$ satisfies
\bequation\label{eq.arc-subtermination}
  \|g_k + (H_k + \sigma_k \|s_k\| I) s_k\| \leq \kappa_3 \|s_k\|^2,
\eequation
then it is used as the trial step.  Otherwise, the process continues with a larger subspace.  We remark that condition \eqref{eq.arc-subtermination} is more restrictive than our condition \eqref{eq.residual_error}, but we use it since it is one that has been proposed for cubic regularization methods; e.g., see (2.13) in \cite{BirgGardMartSantToin17}.

One could employ more sophisticated techniques for setting the elements of the sequence $\{\sigma_k\}$ in \iARC{} that attempt to reduce the number of rejected steps; e.g., see \cite{GoulPorcToin12}.  Such improvements might aid \iRNewton{} as well.  However, for simplicity and to avoid the need for additional parameter tuning, we did not include such enhancements in our implemented algorithms.

As for \iRNewton{}, for consistency between the two implementations, we do not explicitly compute the sequence $\{\lambda_k\}$, but rather employ $\{\sigmal_k\|s_k\|\}$ in its place.  For example, whenever an acceptable step is computed with $\sigmal_k = 0$, then, as described in \textbf{Case 2} in \S\ref{subsec.hybrid}, we effectively use $\lambda_k = 0$.  On the other hand, when $\sigmal_k > 0$, we employ the same iterative approach as used for $\iARC{}$ to compute the trial step~$s_k$ as an approximate minimizer of $c_k(\cdot;\sigmal_k)$, where in place of $\lambda_k$ in \eqref{eq.subproblem_accuracy} we employ $\sigmal_k\|s_k\|$.  Then, in either case, in the remainder of iteration $k \in \N{}_+$, specifically for setting $\sigmal_{k+1}$ and $\sigmau_{k+1}$, we use $\sigmal_k\|s_k\|$ in place of $\lambda_k$ in Steps~\ref{step.reject.lambda} and \ref{step.reject.lambda_large}.  We also define an auxiliary sequence $\{\sigma_k\}$ using the update
\bequationn
  \sigma_{k+1} \gets \bcases \max\left\{\sigmamin,\gamma_0\sigma_k\right\} & \text{if $\rho_k \geq \eta_1$ and $\sigmal_k > 0$} \\[5pt] \sigma_k & \text{if $\sigmal_k = 0$} \\[5pt] \min\left\{\gammal\sigma_k,\sigmamax\right\} & \text{if $\rho_k < \eta_1$ and $\sigmal_k > 0$.} \ecases
\eequationn
This update is similar to the one employed for \iARC{} with the added assurance that $\{\sigma_k\} \subset [\sigmamin,\sigmamax]$.  The elements of this sequence are used in two circumstances.  First, if, as described in \textbf{Case 2} in \S\ref{subsec.hybrid}, CG fails to produce a trial step $s_k$ satisfying \eqref{eq.subproblem_accuracy} (with $\lambda_k = 0$), then we reset $\sigmal_k \gets \sigma_k$ and revert to the same scheme as above to compute the trial step when $\sigmal_k > 0$.  Second, if a step is rejected and $\sigmal_k < \sigmamin$ (equivalently, $\lambda_k < \sigmamin\|s_k\|_2$ as in Step~\ref{step.reject.lambda_small} in Algorithm~\ref{alg.inexact_regularized_Newton}), then we set $\sigmal_{k+1} \gets \sigma_{k+1}$.  Lastly, we note that if CG ever performs $n$ iterations and the resulting solution (due to numerical error) does not satisfy~\eqref{eq.subproblem_accuracy} and no negative curvature is detected, then the resulting approximate solution~$s_k$ is used as the trial step.

\subsection{Results on the \texttt{CUTE}st Test Set}\label{subsec.results}

We employed our implemented algorithms, \iARC{} and \iRNewton{}, to solve unconstrained problems in the \texttt{CUTEst} test set; see \cite{GoulOrbaToin13}.  Among 171 unconstrained problems in the set, one (\texttt{FLETCBV2}) was removed since the algorithms terminated at the initial point, five (\texttt{ARGLINC}, \texttt{DECONVU}, \texttt{FLETCHBV}, \texttt{INDEFM}, and \texttt{POWER}) were removed due to a function evaluation error or our memory limitation of 8GB, and nine (\texttt{EIGENBLS}, \texttt{EIGENCLS}, \texttt{FMINSURF}, \texttt{NONMSQRT}, \texttt{SBRYBND}, \texttt{SCURLY10}, \texttt{SCURLY20}, \texttt{SCURLY30}, and \texttt{SSCOSINE}) were removed since neither algorithm terminated within our time limit.  In addition, four were removed since neither of the algorithms terminated successfully: for \texttt{HIELOW}, \iARC{} reached our maximum iteration limit; for \texttt{CURLY20} and \texttt{SCOSINE}, \iARC{} reached the time limit; for \texttt{INDEF}, \iARC{} terminated due to a subproblem solver error; and for all of these four problems, \iRNewton{} terminated due to our minimum step norm limit.  The remaining set consisted of 152 test problems with number of variables ranging from $2$ to $100,\!000$.  For additional details on the problems used and their sizes, see Appendix~\ref{app.numerical}.


To compare the performance of the implemented algorithms, we generated performance profiles for the number of iterations and number of Hessian-vector products required before termination.  These are shown in Figure~\ref{fig.pp}.  A performance profile graph of an algorithm at point $\alpha$ shows the fraction of the test set for which the algorithm is able to solve within a factor of $2^\alpha$ of the best algorithm for the given measure; see \cite{DolaMore02}.  When generating the profiles, we did not include three of the test problems---\texttt{CURLY10}, \texttt{CURLY30}, and \texttt{MODBEALE}---on which \iARC{} was unsuccessful while \iRNewton{} was successful.  (In particular, \iARC{} reached the time limit for all problems.)  We feel that this gives a fairer comparison with respect to the problems on which both algorithms were successful.

As seen in Figure~\ref{fig.pp}, the algorithms performed relatively comparably when it came to the number of iterations required, though clearly \iRNewton{} had an edge in terms of requiring fewer iterations on various problems.  The difference in terms of numbers of Hessian-vector products required was more drastic, and indeed we point to this as the main measure of improved performance for \iRNewton{} versus~\iARC{}.  One reason for this discrepancy is that \iRNewton{} required fewer iterations on some problems.  However, more significantly, the difference was due in part to \iRNewton's ability to employ and accept inexact Newton steps (with $\lambda_k=0$) on many iterations.  This is due to the fact that, in CG, one is able to compute the Hessian-vector product $H_ks_k$, needed to check the termination conditions for the computation of $s_k$, by taking a linear combination of Hessian-vector products already computed in CG; i.e., if $\{p_{k,i}\}$ are the search directions computed in CG such that $s_k = \sum_i \alpha_{k,i}p_{k,i}$, then CG involves computing $H_kp_{k,i}$ for each $i$ and can compute $H_ks_k = \sum_i \alpha_{k,i} (H_kp_{k,i})$.  By contrast, one is unable to retrieve this product via a linear combination when the step is computed from the minimization of a cubic function, as is needed in \iARC{} and in \iRNewton{} whenever $\sigmal_k > 0$.  Overall, we claim that the primary strength of \iRNewton{} as compared to \iARC{} is its ability to employ inexact Newton steps.

For further details of our numerical results, see Appendix~\ref{app.numerical}.  In these results, we also indicate the number of tridiagonal factorizations required; at least one is needed involving a tridiagonal matrix of size $m \times m$ every time an algorithm solves a cubic subproblem over an $m$-dimensional subspace.

\bfigure
  \centering
  \includegraphics[width=0.38\textwidth,clip=true,trim=35 20 55 15]{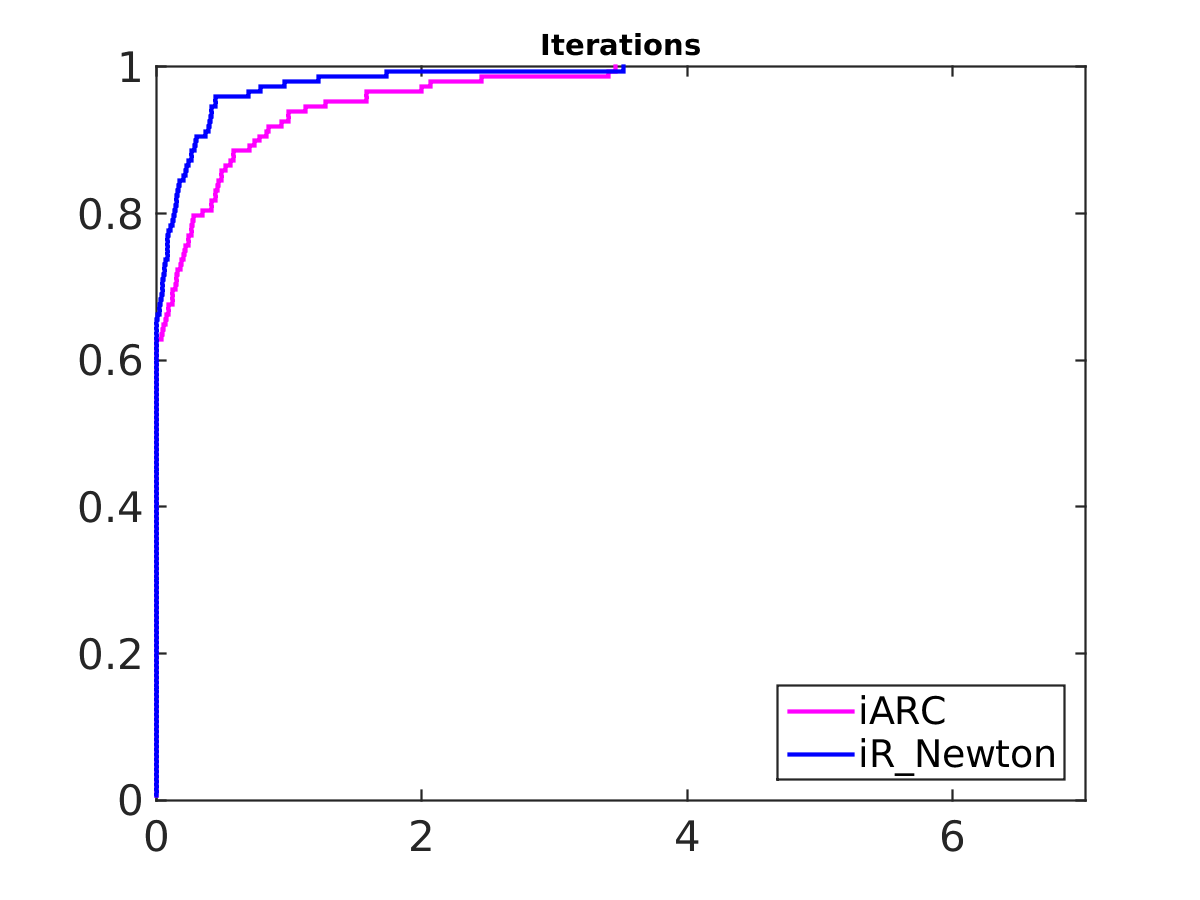}\qquad 
    \includegraphics[width=0.38\textwidth,clip=true,trim=35 20 55 15]{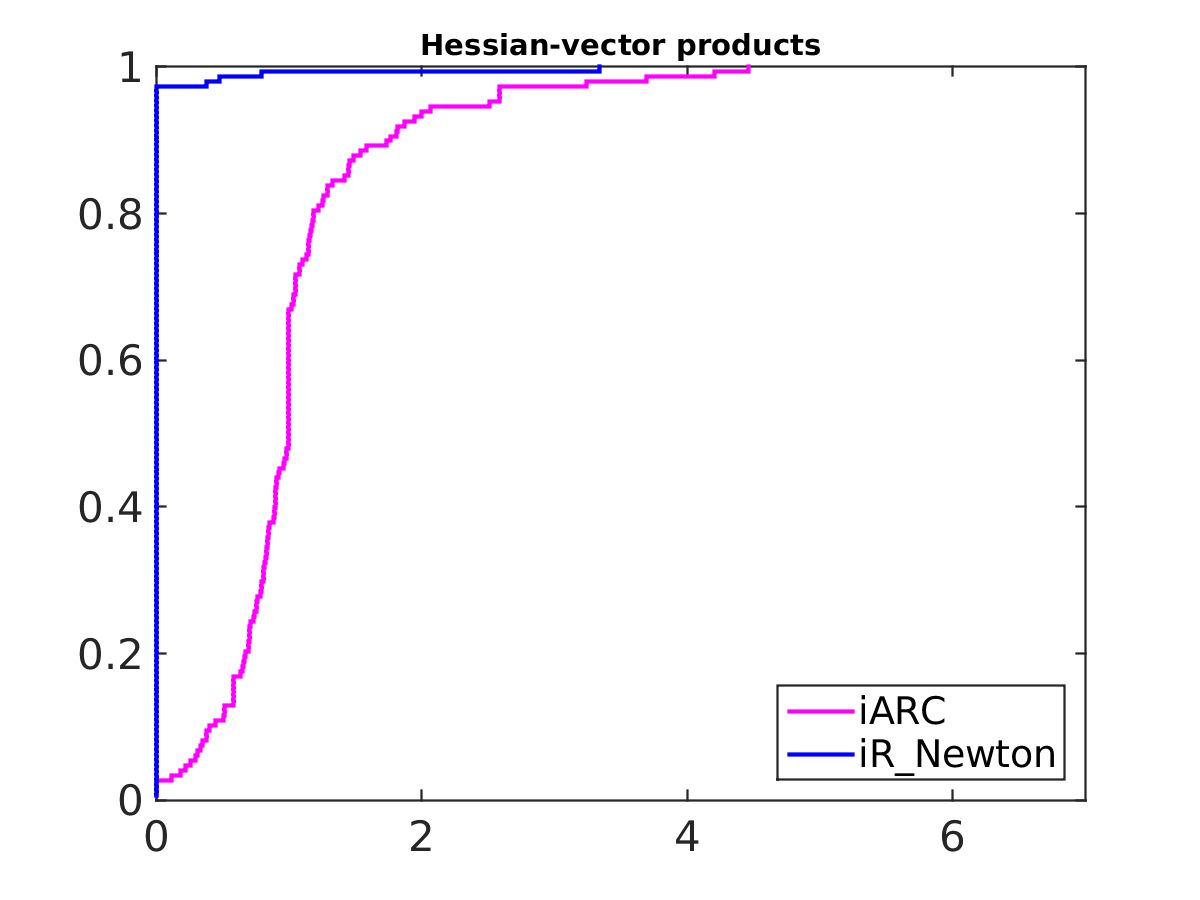}
  \caption{Performance profiles for \iARC{} and \iRNewton.}
  \label{fig.pp}
\efigure

\section{Conclusion}\label{sec.conclusion}

We have proposed a general framework for solving smooth nonconvex optimization problems and proceeded to prove worst-case iteration complexity bounds for it.  In fact, for a certain class of second-order methods employed to minimize a certain class of nonconvex functions, our first-order complexity result for our method is known to be optimal; see~\cite{CartGoulToin11c}.  Our framework is flexible enough to cover a wide range of popular algorithms, an achievement made possible by the use of generic conditions that each trial step is required to satisfy.  The use of such conditions allows for the calculation of inexact Newton steps, for example by performing minimization over expanding Krylov subspaces.  Although we have presented a particular instance of our framework motivated by subproblem~\eqref{prob.primal-dual}, additional instances can easily be derived by applying other optimization strategies for solving \eqref{prob.primal-dual}.  Numerical experiments with an instance of our algorithm showed that it can lead to improved performance on a broad test set as compared to an implementation of a straightforward cubic regularization approach.

\section*{Acknowledgements}

We thank the anonymous referees for their valuable comments, suggestions, and corrections which helped to improve the paper.  We are also grateful to the Associate Editor for handling the paper.

\appendix
\numberwithin{equation}{section}
\numberwithin{theorem}{section}
\numberwithin{lemma}{section}

\section{Subproblem Solution Properties}\label{app.subproblem}

In this appendix, we explore properties of any first-order stationary solution of problem $\Pcal_k(\sigmal_k,\sigmau_k)$ defined as~\eqref{prob.primal-dual}.  Let us define a Lagrangian function for \eqref{prob.primal-dual} as
\bequalin
  \Lcal(s,\lambda,\betal,\betau,\betan) &= f_k + g_k^T s + \thalf s^T (H_k+\lambda I)s \\ &\qquad - \tfrac{\betal}{2}(\lambda^2-(\sigmal_k)^2 \|s\|^2) + \tfrac{\betau}{2}(\lambda^2-(\sigmau_k)^2 \|s\|^2) - \betan \lambda,
\eequalin
where $(\betal,\betau) \in \R{}_+ \times \R{}_+$ are the dual variables associated with the left-hand and right-hand constraints on $\lambda$, respectively, and $\betan \in \R{}_+$ is the dual variable associated with the nonnegativity constraint on~$\lambda$.  The tuple $(s_k,\lambda_k,\betal_k,\betau_k,\betan_k)$ is a first-order primal-dual stationary solution of $\Pcal_k(\sigmal_k,\sigmau_k)$ if it satisfies the following conditions:
\bsubequations\label{eq.primal-dual.kkt}
  \begin{align}
    g_k + (H_k + \lambda_k I) s_k + \betal_k (\sigmal_k)^2 s_k - \betau_k (\sigmau_k)^2 s_k&= 0, \label{eq.primal-dual.kkt.dual1} \\
    \thalf \|s_k\|^2 - \lambda_k (\betal_k - \betau_k) - \betan_k &= 0, \label{eq.primal-dual.kkt.dual2}\\
    0 \leq \betal_k \perp (\lambda_k^2 - (\sigmal_k)^2 \|s_k\|^2) &\geq 0, \label{eq.primal-dual.kkt.comp1} \\
    0 \leq \betau_k \perp (\lambda_k^2 - (\sigmau_k)^2 \|s_k\|^2) &\leq 0, \ \ \text{and}  \label{eq.primal-dual.kkt.comp2} \\
    0 \leq \betan_k \perp \lambda_k & \geq 0. \label{eq.primal-dual.kkt.comp3}
  \end{align}
\esubequations

We make the following assumption throughout this appendix.

\bassumption\label{ass.g-nonzero}
  \textit{
  The vector $g_k$ is nonzero.
  }
\eassumption

Under this assumption, the following lemma is a simple consequence of~\eqref{eq.primal-dual.kkt.dual1}.

\blemma\label{lem.s-nonzero}
  \textit{
  Any solution of \eqref{prob.primal-dual} has $s_k \neq 0$.
  }
\elemma

We now establish conditions that must hold depending on the value of $\sigmal_k \in \R{}_+$.

\blemma\label{lem.sigmal}
  \textit{
  The following hold true for any solution of \eqref{eq.primal-dual.kkt}.
  \bitemize
    \item[(i)] If $\sigmal_k > 0$, then $\lambda_k > 0$, $\betan_k = 0$, $\betal_k > 0$, and $\lambda_k = \sigmal_k \|s_k\|$.
    \item[(ii)] If $\sigmal_k = 0$, then $\lambda_k = 0$.
  \eitemize
  }
\elemma
\bproof
  Consider part $(i)$.  For the sake of deriving a contradiction, suppose $\sigmal_k > 0$ and $\lambda_k = 0$.  These, along with Lemma~\ref{lem.s-nonzero}, imply that $0 = \lambda_k^2 < (\sigmal_k)^2 \|s_k\|^2$, which contradicts~\eqref{eq.primal-dual.kkt.comp1}.  Hence, $\lambda_k > 0$, as claimed.  Then, it follows from \eqref{eq.primal-dual.kkt.comp3} that $\betan_k = 0$, as claimed.  Next, observe that from~\eqref{eq.primal-dual.kkt.dual2}, Lemma~\ref{lem.s-nonzero}, $\betan_k = 0$, $\lambda_k > 0$, and $(\betal_k,\betau_k) \geq 0$, it follows that $\betal_k > 0$, as claimed.  This, along with~\eqref{eq.primal-dual.kkt.comp1}, implies that $\lambda_k^2 = (\sigmal_k)^2 \|s_k\|^2$. This implies that $\lambda_k = \pm (\sigmal_k)\|s_k\|$, which combined with $\lambda_k \in\Re^{+}$ means that $\lambda_k = \sigmal_k \|s_k\|$, as claimed.

  Now consider part $(ii)$.  For the sake of deriving a contradiction, suppose that $\sigmal_k = 0$ and $\lambda_k > 0$. Then, it follows from \eqref{eq.primal-dual.kkt.comp3} that $\betan_k = 0$.  Moreover, combining $\sigmal_k = 0$ and $\lambda_k > 0$, it follows from \eqref{eq.primal-dual.kkt.comp1} that $\betal_k = 0$.  It now follows from $\betal_k = 0$, $\betan_k = 0$, and \eqref{eq.primal-dual.kkt.dual2} that
\bequation \label{eq.s-beta2-lambda}
  \thalf\|s_k\|^2 = - \lambda_k \betau_k \leq 0,
\eequation
where the inequality follows from $\lambda_k > 0$ and $\betau_k \geq 0$.  This contradicts Lemma~\ref{lem.s-nonzero}.
\eproof

Our main result is the following.  In part~$(i)$ with $\sigmal_k > 0$, we show that solving \eqref{prob.primal-dual} is equivalent to solving what may be referred to as an \ARC\ subproblem~\cite{CartGoulToin11a}.  In part~$(ii)$ with $\sigmal_k = 0$, we show that it is equivalent to minimizing a quadratic, if a minimizer exists.

\btheorem\label{th.equivalency}
  \textit{
  The following hold true.
  \bitemize
    \item[(i)] Suppose $\sigmal_k > 0$.  Then, \eqref{prob.primal-dual} has a solution $(s_k,\lambda_k)$, which can be obtained as
    \bequation\label{prob.primal-dual.sigma-reduced}
      s_k \in \arg\min_{s \in \R{n}}\ (f_k + g_k^T s + \thalf s^T H_k s + \thalf \sigmal_k \|s\|^3),
    \eequation
    then setting $\lambda_k = \sigmal_k \|s_k\| > 0$.
    \item[(ii)] If $\sigmal_k = 0$, then a solution of problem~\eqref{prob.primal-dual} exists if and only if $H_k \succeq 0$ and $g_k^Tu = 0$ for all $u \in \Null(H_k)$.  In such cases, computing a solution $(s_k,\lambda_k)$ of problem~\eqref{prob.primal-dual} is equivalent to computing a solution $s_k$ of problem~\eqref{prob.quadratic} and setting $\lambda_k = 0$.
  \eitemize
  }
\etheorem
\bproof
  Consider part $(i)$.  Since $\sigmal_k > 0$, it follows from Lemma~\ref{lem.sigmal} that problem~\eqref{prob.primal-dual} is equivalent to
  \bequation\label{prob.primal-dual.sigmal}
  \baligned 
    \min_{(s,\lambda) \in \R{n}\times\R{}_{+}}~&f_k + g_k^T s + \thalf s^T (H_k + \lambda I) s \\
    \st~~ & \sigmal_k \|s\| = \lambda,
  \ealigned
  \eequation
  where, by Lemma~\ref{lem.s-nonzero}, it follows that the solution has $\lambda_k > 0$, as desired.  Substituting the constraint of~\eqref{prob.primal-dual.sigmal} into the objective of \eqref{prob.primal-dual.sigmal}, one finds that solving it is equivalent to solving~\eqref{prob.primal-dual.sigma-reduced} for $s_k$, then setting $\lambda_k = \sigmal_k \|s_k\|$, as claimed.  Since $\sigmal_k > 0$, a minimizer of problem~\eqref{prob.primal-dual.sigma-reduced} exists because it involves the minimization of a coercive function.

  Now consider part $(ii)$.  Since $\sigmal_k = 0$, it follows from Lemma~\ref{lem.sigmal} that $\lambda_k = 0$, meaning that problem~\eqref{prob.primal-dual} is equivalent to \eqref{prob.quadratic}.  This problem has a solution if and only if the objective is bounded below, which is the case if and only if $H_k \succeq 0$ and $g_k^Tu = 0$ for all $u \in \Null(H_k)$.
\eproof

\section{Subproblem Solution Properties Over Subspaces}\label{app.subproblem.reduced}

In this appendix, we explore properties of any first-order stationary solution (when one exists) of problem $\Pcal_k(\sigmal_k,\sigmau_k)$ defined as \eqref{prob.primal-dual} when the search space for $s$ is restricted to a subspace of~$\R{n}$.  Specifically, for some $m$-dimensional subspace $\Lcal_k \subseteq \R{n}$, consider the problem
\bequation\label{prob.primal-dual_subspace}
  \baligned
    \min_{(s,\lambda) \in \Lcal_k \times\R{}_{+}} &\ f_k + g_k^T s + \thalf s^T (H_k + \lambda I) s \\
    \st &\ (\sigmal_k)^2 \|s\|^2 \leq \lambda^2 \leq (\sigmau_k)^2 \|s\|^2.
  \ealigned
\eequation
Given an orthogonal basis $R_k$ for $\Lcal_k$, a solution of~\eqref{prob.primal-dual_subspace} can be obtained from that of
\bequation\label{prob.primal-dual_subspace-reduced}
  \baligned
    \min_{(v,\lambda) \in \R{m} \times\R{}_{+}} &\ f_k + g_k^T R_k v + \thalf (R_kv)^T (H_k + \lambda I) R_kv \\
    \st &\ (\sigmal_k)^2 \|v\|^2 \leq \lambda^2 \leq (\sigmau_k)^2 \|v\|^2.
  \ealigned
\eequation
Specifically, if $(v_k,\lambda_k,\betal_k,\betau_k,\betan_k)$ is a first-order primal-dual solution of problem~\eqref{prob.primal-dual_subspace-reduced}, then the tuple $(s_k,\lambda_k,\betal_k,\betau_k,\betan_k)$ with $s_k = R_kv_k$ is such a  solution of problem~\eqref{prob.primal-dual_subspace}.

In Appendix~\ref{app.subproblem}, we proved properties of a solution (if one exists) of a problem of the form \eqref{prob.primal-dual_subspace-reduced}.  Let us now translate the results of that appendix to the present setting, for which we require the following assumption on the reduced gradient $R_k^Tg_k$.

\bassumption\label{ass.Zg-nonzero}
  \textit{
  The vector $R_k^Tg_k$ is nonzero.
  }
\eassumption

\blemma\label{lem.v-nonzero}
  \textit{
  Any solution of \eqref{prob.primal-dual_subspace-reduced} has $v_k \neq 0$.
  }
\elemma

\blemma\label{lem.sigmal-v}
  \textit{
  The following hold for any first-order primal-dual solution of~\eqref{prob.primal-dual_subspace}.
  \bitemize
    \item[(i)] If $\sigmal_k > 0$, then $\lambda_k > 0$, $\betan_k = 0$, $\betal_k > 0$, and $\lambda_k = \sigmal_k \|v_k\|$.
    \item[(ii)] If $\sigmal_k = 0$, then $\lambda_k = 0$.
  \eitemize
  }
\elemma

\btheorem\label{th.equivalency_subspace}
  \textit{
  The following hold true.
  \bitemize
    \item[(i)] Suppose $\sigmal_k > 0$.  Then, \eqref{prob.primal-dual_subspace-reduced} has a solution $(v_k,\lambda_k)$, which can be obtained as
    \bequation\label{prob.primal-dual_subspace.sigma-reduced}
      v_k \in \arg\min_{v \in \R{m}}\ (f_k + g_k^T R_k v + \thalf v^T R_k^T H_k R_k v + \thalf \sigmal_k \|v\|^3),
    \eequation
    then setting $\lambda_k = \sigmal_k \|v_k\| > 0$.
    \item[(ii)] If $\sigmal_k = 0$, then a solution of~\eqref{prob.primal-dual_subspace-reduced} exists if and only if $R_k^TH_kR_k \succeq 0$ and $g_k^TR_ku = 0$ for all $u \in \Null(R_k^T H_k R_k)$.  In such cases, computing a solution $(v_k,\lambda_k)$ of problem~\eqref{prob.primal-dual_subspace-reduced} is equivalent to computing a solution $v_k$ of
    \bequation\label{prob.v}
      \min_{v\in\R{m}}\ f_k + g_k^TR_kv + \thalf v^T R_k^TH_kR_k v
    \eequation
    and setting $\lambda_k = 0$.
  \eitemize
  }
\etheorem

Considering problem~\eqref{prob.primal-dual_subspace.sigma-reduced}, we obtain the following result from \cite[Lemma 3.2]{CartGoulToin11a}.

\blemma\label{lem.cubicsolve}
  \textit{
  If $\sigmal_k > 0$, then $v_k$ from \eqref{prob.primal-dual_subspace.sigma-reduced} satisfies
  \bsubequations
    \begin{align}
      g_k^TR_kv_k + v_k^T R_k^T H_k R_k v_k + \tfrac32 \sigmal_k \|v_k\|^3 &= 0 \label{eq.cubicsolve1} \\
      v_k^T R_k^T H_k R_k v_k + \tfrac32 \sigmal_k \|v_k\|^3 &\geq 0 \label{eq.cubicsolve2} \\
      R_k^T H_k R_k + \tfrac32 \sigmal_k \|v_k\| I &\succeq 0. \label{eq.cubicsolve3}
    \end{align}
  \esubequations
  }
\elemma

We now show that, under certain reasonable assumptions, solutions of the primal-dual reduced-space subproblem~\eqref{prob.primal-dual_subspace} satisfy the conditions required by Assumptions~\ref{ass.subproblem_accuracy} and \ref{ass.2nd}. 

\btheorem\label{th.consistency-subspace}
  \textit{
  The following hold true.
  \bitemize
    \item[(a)] Any solution of problem~\eqref{prob.primal-dual_subspace} satisfies \eqref{eq.subspace_optimality}.  
    \item[(b)] Any solution of problem~\eqref{prob.primal-dual_subspace} satisfies \eqref{eq.Cauchy_decrease} provided $g_k \in \Lcal_k$.
    \item[(c)] Any solution of problem~\eqref{prob.primal-dual_subspace} satisfies \eqref{eq.residual_error} provided $\Lcal_k = \R{n}$.
    \item[(d)] Any solution of problem~\eqref{prob.primal-dual_subspace} satisfies \eqref{eq.2nd} for any $\kappa_4 \geq \tfrac32\sup_{k\in\N{}_+}\{\sigmal_k\}$.
  \eitemize
  }
\etheorem
\bproof
  Any first-order primal-dual solution $(s_k,\lambda_k,\betal_k,\betau_k,\betan_k)$ of problem~\eqref{prob.primal-dual_subspace} corresponds to such a solution $(v_k,\lambda_k,\betal_k,\betau_k,\betan_k)$ of problem~\eqref{prob.primal-dual_subspace-reduced} where $s_k = R_kv_k$.  Hence, throughout this proof, for any solution vector $s_k$ for problem~\eqref{prob.primal-dual_subspace}, we may let $s_k = R_kv_k$ where $v_k$ satisfies the properties in Lemmas~\ref{lem.v-nonzero}--\ref{lem.cubicsolve}.

  First, suppose $\sigmal_k > 0$, which by Theorem~\ref{th.equivalency_subspace}$(i)$ implies that problem~\eqref{prob.primal-dual_subspace} has a solution.  Then, it follows from \eqref{eq.cubicsolve1}, $s_k = R_k v_k$, and Lemma~\ref{lem.sigmal-v}(i) that
  \bequationn
    0 = g_k^T s_k + s_k^T H_k s_k + \tfrac{3}{2} \sigmal_k \|s_k\|^3 = g_k^T s_k + s_k^T H_k s_k + \tfrac{3}{2} \lambda_k \|s_k\|^2,
  \eequationn
  which means that
  \bequation\label{eq.like-decrease}
    s_k^T (g_k + (H_k+\lambda_k I) s_k) = -\thalf \lambda_k \|s_k\|^2.
  \eequation
  Meanwhile, from \eqref{eq.cubicsolve2}, $s_k = R_k v_k$, and~Lemma~\ref{lem.sigmal-v}$(i)$, it follows that
  \bequationn
    0 \leq s_k^T H_k s_k + \tfrac{3}{2} \sigmal_k \|s_k\|^3 = s_k^T H_k s_k + \tfrac{3}{2} \lambda_k \|s_k\|^2 = s_k^T (H_k + \lambda_k I)s_k + \thalf \lambda_k \|s_k\|^2,
  \eequationn
  which means that
  \bequation  \label{eq.like-psd}
   -\tfrac{1}{4} \lambda_k \|s_k\|^2 
   \leq \thalf s_k^T (H_k + \lambda_k I)s_k.
  \eequation
  It follows from~\eqref{eq.like-decrease}, \eqref{eq.like-psd}, $\lambda_k > 0$ (by Lemma~\ref{lem.sigmal-v}$(i)$), and $(\kappa_1,\kappa_2) \in \R{}_{++} \times \R{}_{++}$ that
  \begin{align*}
    s_k^T (g_k + (H_k+\lambda_k I) s_k) = -\thalf \lambda_k \|s_k\|^2
      &\leq \min\{\thalf \kappa_1\|s_k\|^2, \thalf s_k^T (H_k + \lambda_k I)s_k -\tfrac{1}{4} \lambda_k \|s_k\|^2 \} \\
      &\leq \min\{\thalf \kappa_1\|s_k\|^2, \thalf s_k^T (H_k + \lambda_k I)s_k + \thalf \kappa_2\|s\|^3 \},
  \end{align*}
  which implies \eqref{eq.subspace_optimality}.  This establishes that part~$(a)$ is true.  Now consider part $(b)$.  From Theorem~\ref{th.equivalency_subspace}, \cite[Lemma 2.1]{CartGoulToin11a}, and $s_k = R_k v_k$, it follows that 
  \bequationn
    f_k - q_k(s_k) - \thalf \sigmal_k \|s_k\|^3
    \geq \frac{\|R_k^T g_k\|}{6\sqrt{2}} \min \left\{\frac{\|R_k^Tg_k\|}{1+\|R_k^TH_kR_k\|}, \frac{1}{\sqrt{6}}\sqrt{\frac{\|R_k^Tg_k\|}{\sigmal_k}}\right\}.
  \eequationn
  Since, under assumption, $g_k \in \Lcal_k$ so that $g_k = R_k y$ for some $y\in\R{m}$, it follows that
  \bequationn
    \|R_k^T g_k\| = \|R_k^T R_k y\| = \|y\| = \|R_ky\| = \|g_k\|.
  \eequationn
  Combining this with $\|R_k^T H_k R_k\| \leq \|H_k\|$ and the previous displayed inequality shows
  \bequationn
    f_k - q_k(s_k) - \thalf \sigmal_k \|s_k\|^3
    \geq \frac{\|g_k\|}{6\sqrt{2}} \min \left\{\frac{\|g_k\|}{1+\|H_k\|}, \frac{1}{\sqrt{6}}\sqrt{\frac{\|g_k\|}{\sigmal_k}}\right\}.
  \eequationn
This may now be combined with Theorem~\ref{th.equivalency_subspace} (specifically $\lambda_k = \sigmal_k \|s_k\| > 0$) to obtain
  \bequationn
    f_k - q_k(s_k) 
    \geq f_k - q_k(s_k) - \thalf \sigmal_k \|s_k\|^3
    \geq \frac{\|g_k\|}{6\sqrt{2}} \min \left\{\frac{\|g_k\|}{1+\|H_k\|}, \frac{1}{\sqrt{6}}\sqrt{\frac{\|g_k\|\|s_k\|}{\lambda_k}}\right\},
  \eequationn
  which means that $(s_k,\lambda_k)$ satisfies \eqref{eq.Cauchy_decrease}, proving part $(b)$.  Now consider part $(c)$.  It follows from Theorem~\ref{th.equivalency}$(i)$ and the optimality conditions for problem~\eqref{prob.primal-dual.sigma-reduced} that
  \bequationn
    0 = g_k + H_ks_k + \tfrac32\sigmal_k \|s_k\| s_k = g_k + H_ks_k + \tfrac32 \lambda_k s_k = g_k + (H_k + \lambda_k I)s_k + \thalf \lambda_k s_k.
  \eequationn
  This and the fact that $\kappa_3 > 0$ imply that
  \bequationn
    \baligned
      \|g_k + (H_k + \lambda_k I)s_k\| = \thalf \lambda_k \|s_k\| \leq \lambda_k\|s_k\| + \kappa_3\|s_k\|^2,
    \ealigned
  \eequationn
  which completes the proof of part $(c)$.  Finally, consider part $(d)$.  From \eqref{eq.cubicsolve3}, the fact that $\|s_k\| = \|v_k\|$, and $\kappa_4 \geq \tfrac32\sup_{k\in\N{}_+}\{\sigmal_k\}$, it follows that
  \bequationn
    \xi(R_k^TH_kR_k) \geq -\tfrac32\sigmal_k\|s_k\| \geq -\kappa_4\|s_k\|,
  \eequationn
  as desired to prove part $(d)$.
  
  Now suppose that $\sigmal_k = 0$.  From Theorem~\ref{th.equivalency_subspace}$(ii)$, a solution of problem~\eqref{prob.primal-dual_subspace} exists if and only if $R_k^TH_kR_k \succeq 0$ and $g_k^TR_ku = 0$ for all $u \in \Null(R_k^TH_kR_k)$.  If this is not the case, then there is nothing left to prove; hence, let us assume that these conditions hold.  From these conditions, Theorem~\ref{th.equivalency_subspace}$(ii)$, the optimality conditions of problem~\eqref{prob.v}, the fact that $\lambda_k = 0$, and $s_k = R_kv_k$, it follows that
  \bequationn
    g_k^Ts_k + s_k^TH_ks_k = 0\ \ \text{and}\ \ s_k^TH_ks_k \geq 0.
  \eequationn
  This shows that \eqref{eq.subspace_optimality} holds, proving part $(a)$ for this case.  Next, since $v_k$ is given by the solution of problem~\eqref{prob.v}, it follows that the reduction in the objective yielded by $v_k$ is at least as large as the reduction obtained by minimizing the objective over the span of $-R_k^Tg_k$.  Hence, from standard theory on Cauchy decrease (see \cite{ConnGoulToin00} or \cite{NoceWrig06}), one can conclude that
  \bequationn
    f_k - q_k(s_k) \geq \frac{\|R_k^Tg_k\|}{2} \min\left\{ \frac{\|R_k^Tg_k\|}{1 + \|R_k^TH_kR_k\|} , \|s_k\| \right\}.
  \eequationn
  Thus, using the arguments in the previous paragraph under the assumption that $g_k \in \Lcal_k$, one is led to the conclusion that \eqref{eq.Cauchy_decrease} holds, which proves part $(b)$ for this case.  Next, when $\Lcal_k = \R{n}$, the optimality conditions for problem~\eqref{prob.v} imply that $g_k + H_k s_k = 0$, which, since $\lambda_k = 0$, implies that \eqref{eq.residual_error} holds, proving part $(c)$.  Finally, since $R_k^TH_kR_k \succeq 0$, it follows that \eqref{eq.2nd} holds, proving part $(d)$.
\eproof

\section{Detailed Numerical Results}\label{app.numerical}

Further details of the results of our numerical experiments are shown in Table~\ref{tbl.numericalresults}. In the table, \texttt{\#Var} indicates the number of variables, \texttt{\#Iter} indicates the number of iterations required (with \texttt{\%Newton} indicating the percentage that were inexact Newton steps with $\lambda_k=0$), \texttt{\#Acc} indicates the number of accepted steps (again with \texttt{\%Newton} indicating the percentage that were inexact Newton steps), \texttt{\#Hv-prod} indicates the number of Hessian-vector products required, and \texttt{\#T-fact} indicates the number of tridiagonal matrix factorizations required.

\begin{scriptsize}
\renewcommand{\tabcolsep}{0.18cm}
\texttt{
\begin{longtable}{|c|c|l|rr|rr|r|r|}
\caption{Numerical results for \iARC{} and \iRNewton.} \label{tbl.numericalresults} \\
\hline \multicolumn{1}{|c|}{\textbf{Prob}} & \multicolumn{1}{|c|}{\textbf{\#Var}} & \multicolumn{1}{|c|}{\textbf{Alg}} & \multicolumn{2}{|c|}{\textbf{\#Iter} (\textbf{\%Newton})} & \multicolumn{2}{|c|}{\textbf{\#Acc} (\textbf{\%Newton})} & \multicolumn{1}{|c|}{\textbf{\#Hv-prod}} & \multicolumn{1}{|c|}{\textbf{\#T-fact}} \\ \hline
\endfirsthead
\multicolumn{9}{c}{{\bfseries \tablename\ \thetable{} -- continued from previous page}} \\
\hline \multicolumn{1}{|c|}{\textbf{Prob}} & \multicolumn{1}{|c|}{\textbf{\#Var}} & \multicolumn{1}{|c|}{\textbf{Alg}} & \multicolumn{2}{|c|}{\textbf{\#Iter} (\textbf{\%Newton})} & \multicolumn{2}{|c|}{\textbf{\#Acc} (\textbf{\%Newton})} & \multicolumn{1}{|c|}{\textbf{\#Hv-prod}} & \multicolumn{1}{|c|}{\textbf{\#T-fact}} \\ \hline
\endhead
\endfoot
\endlastfoot
    \multirow{ 2}{*}{     AKIVA} &  \multirow{ 2}{*}{    2} &       iARC &       5 &&           5 &&        15 &       20    \\ 
&&  iR\_Newton &       5 &(\%100) &           5 &(\%100) &        10 &        0    \\ \hline 
\multirow{ 2}{*}{  ALLINITU} &  \multirow{ 2}{*}{    4} &       iARC &      11 &&           8 &&        56 &       61    \\ 
&&  iR\_Newton &       8 &(\%50) &           6 &(\%67) &        25 &       21    \\ \hline 
\multirow{ 2}{*}{   ARGLINA} &  \multirow{ 2}{*}{  200} &       iARC &       3 &&           3 &&         6 &        3    \\ 
&&  iR\_Newton &       1 &(\%100) &           1 &(\%100) &         1 &        0    \\ \hline 
\multirow{ 2}{*}{   ARGLINB} &  \multirow{ 2}{*}{  200} &       iARC &       2 &&           2 &&         4 &        2    \\ 
&&  iR\_Newton &       1 &(\%100) &           1 &(\%100) &         1 &        0    \\ \hline 
\multirow{ 2}{*}{   ARWHEAD} &  \multirow{ 2}{*}{ 5000} &       iARC &       4 &&           4 &&        10 &        6    \\ 
&&  iR\_Newton &       4 &(\%100) &           4 &(\%100) &         5 &        0    \\ \hline 
\multirow{ 2}{*}{      BARD} &  \multirow{ 2}{*}{    3} &       iARC &      11 &&           8 &&        50 &       52    \\ 
&&  iR\_Newton &      11 &(\%91) &          10 &(\%90) &        28 &        6    \\ \hline 
\multirow{ 2}{*}{   BDQRTIC} &  \multirow{ 2}{*}{ 5000} &       iARC &       9 &&           9 &&        34 &       33    \\ 
&&  iR\_Newton &       9 &(\%100) &           9 &(\%100) &        17 &        0    \\ \hline 
\multirow{ 2}{*}{     BEALE} &  \multirow{ 2}{*}{    2} &       iARC &      11 &&           8 &&        33 &       47    \\ 
&&  iR\_Newton &      12 &(\%42) &           8 &(\%62) &        29 &       29    \\ \hline 
\multirow{ 2}{*}{    BIGGS6} &  \multirow{ 2}{*}{    6} &       iARC &     457 &&         383 &&      3446 &     3679    \\ 
&&  iR\_Newton &     419 &(\%96) &         406 &(\%97) &      1557 &      183    \\ \hline 
\multirow{ 2}{*}{       BOX} &  \multirow{ 2}{*}{10000} &       iARC &       3 &&           3 &&        14 &       12    \\ 
&&  iR\_Newton &       4 &(\%75) &           3 &(\%67) &         9 &        4    \\ \hline 
\multirow{ 2}{*}{      BOX3} &  \multirow{ 2}{*}{    3} &       iARC &       7 &&           7 &&        30 &       32    \\ 
&&  iR\_Newton &       7 &(\%100) &           7 &(\%100) &        16 &        0    \\ \hline 
\multirow{ 2}{*}{  BOXPOWER} &  \multirow{ 2}{*}{20000} &       iARC &       3 &&           3 &&        10 &        9    \\ 
&&  iR\_Newton &       7 &(\%100) &           7 &(\%100) &        13 &        0    \\ \hline 
\multirow{ 2}{*}{    BRKMCC} &  \multirow{ 2}{*}{    2} &       iARC &       2 &&           2 &&         6 &        7    \\ 
&&  iR\_Newton &       2 &(\%100) &           2 &(\%100) &         4 &        0    \\ \hline 
\multirow{ 2}{*}{   BROWNAL} &  \multirow{ 2}{*}{  200} &       iARC &       2 &&           2 &&         6 &        4    \\ 
&&  iR\_Newton &       1 &(\%100) &           1 &(\%100) &         1 &        0    \\ \hline 
\multirow{ 2}{*}{   BROWNBS} &  \multirow{ 2}{*}{    2} &       iARC &      53 &&          38 &&       142 &      191    \\ 
&&  iR\_Newton &       5 &(\%80) &           5 &(\%80) &        11 &        5    \\ \hline 
\multirow{ 2}{*}{  BROWNDEN} &  \multirow{ 2}{*}{    4} &       iARC &       8 &&           8 &&        35 &       35    \\ 
&&  iR\_Newton &       9 &(\%100) &           9 &(\%100) &        20 &        0    \\ \hline 
\multirow{ 2}{*}{  BROYDN7D} &  \multirow{ 2}{*}{ 5000} &       iARC &     472 &&         279 &&      7202 &    12598    \\ 
&&  iR\_Newton &     812 &(\%29) &         346 &(\%2) &      5033 &    10022    \\ \hline 
\multirow{ 2}{*}{    BRYBND} &  \multirow{ 2}{*}{ 5000} &       iARC &      19 &&          10 &&       240 &      367    \\ 
&&  iR\_Newton &      17 &(\%53) &          11 &(\%82) &       206 &       81    \\ \hline 
\multirow{ 2}{*}{  CHAINWOO} &  \multirow{ 2}{*}{ 4000} &       iARC &      81 &&          57 &&       798 &      942    \\ 
&&  iR\_Newton &      70 &(\%81) &          59 &(\%88) &       409 &      185    \\ \hline 
\multirow{ 2}{*}{  CHNROSNB} &  \multirow{ 2}{*}{   50} &       iARC &      64 &&          40 &&      1126 &     1456    \\ 
&&  iR\_Newton &      53 &(\%75) &          40 &(\%82) &       499 &      294    \\ \hline 
\multirow{ 2}{*}{  CHNRSNBM} &  \multirow{ 2}{*}{   50} &       iARC &      96 &&          58 &&      1708 &     2320    \\ 
&&  iR\_Newton &     101 &(\%58) &          59 &(\%59) &       899 &      960    \\ \hline 
\multirow{ 2}{*}{     CLIFF} &  \multirow{ 2}{*}{    2} &       iARC &      14 &&          14 &&        28 &       14    \\ 
&&  iR\_Newton &      14 &(\%100) &          14 &(\%100) &        14 &        0    \\ \hline 
\multirow{ 2}{*}{    COSINE} &  \multirow{ 2}{*}{10000} &       iARC &      12 &&           7 &&       108 &      140    \\ 
&&  iR\_Newton &      11 &(\%55) &           7 &(\%71) &        45 &       29    \\ \hline 
\multirow{ 2}{*}{  CRAGGLVY} &  \multirow{ 2}{*}{ 5000} &       iARC &      30 &&          30 &&       228 &      208    \\ 
&&  iR\_Newton &      31 &(\%100) &          31 &(\%100) &       108 &        0    \\ \hline 
\multirow{ 2}{*}{      CUBE} &  \multirow{ 2}{*}{    2} &       iARC &      42 &&          27 &&       126 &      170    \\ 
&&  iR\_Newton &      35 &(\%69) &          25 &(\%76) &        72 &       51    \\ \hline 
\multirow{ 2}{*}{   CURLY10} &  \multirow{ 2}{*}{10000} &       iARC & ---   &&   --- &&   --- &   ---    \\ 
&&  iR\_Newton &     328 &(\%95) &         318 &(\%97) &    271881 &    19957    \\ \hline 
\multirow{ 2}{*}{   CURLY30} &  \multirow{ 2}{*}{10000} &       iARC &  ---  &&   --- &&   --- &   ---    \\ 
&&  iR\_Newton &      87 &(\%83) &          77 &(\%91) &    125639 &      630    \\ \hline 
\multirow{ 2}{*}{  DENSCHNA} &  \multirow{ 2}{*}{    2} &       iARC &       5 &&           5 &&        15 &       15    \\ 
&&  iR\_Newton &       5 &(\%100) &           5 &(\%100) &        10 &        0    \\ \hline 
\multirow{ 2}{*}{  DENSCHNB} &  \multirow{ 2}{*}{    2} &       iARC &       7 &&           6 &&        20 &       22    \\ 
&&  iR\_Newton &       7 &(\%71) &           5 &(\%80) &        12 &        8    \\ \hline 
\multirow{ 2}{*}{  DENSCHNC} &  \multirow{ 2}{*}{    2} &       iARC &      13 &&           9 &&        38 &       46    \\ 
&&  iR\_Newton &      11 &(\%82) &           9 &(\%89) &        21 &        8    \\ \hline 
\multirow{ 2}{*}{  DENSCHND} &  \multirow{ 2}{*}{    3} &       iARC &      61 &&          57 &&       206 &      172    \\ 
&&  iR\_Newton &      44 &(\%86) &          40 &(\%95) &        82 &       27    \\ \hline 
\multirow{ 2}{*}{  DENSCHNE} &  \multirow{ 2}{*}{    3} &       iARC &      24 &&          15 &&        68 &       62    \\ 
&&  iR\_Newton &      21 &(\%52) &          16 &(\%69) &        41 &       58    \\ \hline 
\multirow{ 2}{*}{  DENSCHNF} &  \multirow{ 2}{*}{    2} &       iARC &       5 &&           5 &&        15 &       15    \\ 
&&  iR\_Newton &       5 &(\%100) &           5 &(\%100) &        10 &        0    \\ \hline 
\multirow{ 2}{*}{  DIXMAANA} &  \multirow{ 2}{*}{ 3000} &       iARC &       6 &&           6 &&        14 &        8    \\ 
&&  iR\_Newton &       6 &(\%100) &           6 &(\%100) &         7 &        0    \\ \hline 
\multirow{ 2}{*}{  DIXMAANB} &  \multirow{ 2}{*}{ 3000} &       iARC &       7 &&           7 &&        16 &        9    \\ 
&&  iR\_Newton &       7 &(\%100) &           7 &(\%100) &         8 &        0    \\ \hline 
\multirow{ 2}{*}{  DIXMAANC} &  \multirow{ 2}{*}{ 3000} &       iARC &       8 &&           8 &&        18 &        9    \\ 
&&  iR\_Newton &       8 &(\%100) &           8 &(\%100) &         9 &        0    \\ \hline 
\multirow{ 2}{*}{  DIXMAAND} &  \multirow{ 2}{*}{ 3000} &       iARC &       9 &&           9 &&        20 &       10    \\ 
&&  iR\_Newton &       9 &(\%100) &           9 &(\%100) &        10 &        0    \\ \hline 
\multirow{ 2}{*}{  DIXMAANE} &  \multirow{ 2}{*}{ 3000} &       iARC &      59 &&          59 &&       670 &      622    \\ 
&&  iR\_Newton &      60 &(\%100) &          60 &(\%100) &       331 &        0    \\ \hline 
\multirow{ 2}{*}{  DIXMAANF} &  \multirow{ 2}{*}{ 3000} &       iARC &      38 &&          37 &&       510 &      487    \\ 
&&  iR\_Newton &      37 &(\%100) &          37 &(\%100) &       249 &        0    \\ \hline 
\multirow{ 2}{*}{  DIXMAANG} &  \multirow{ 2}{*}{ 3000} &       iARC &      39 &&          39 &&       532 &      514    \\ 
&&  iR\_Newton &      40 &(\%100) &          40 &(\%100) &       288 &        0    \\ \hline 
\multirow{ 2}{*}{  DIXMAANH} &  \multirow{ 2}{*}{ 3000} &       iARC &      41 &&          41 &&       448 &      421    \\ 
&&  iR\_Newton &      41 &(\%100) &          41 &(\%100) &       224 &        0    \\ \hline 
\multirow{ 2}{*}{  DIXMAANI} &  \multirow{ 2}{*}{ 3000} &       iARC &     193 &&         193 &&      3456 &     3464    \\ 
&&  iR\_Newton &     249 &(\%100) &         249 &(\%100) &      2814 &        0    \\ \hline 
\multirow{ 2}{*}{  DIXMAANJ} &  \multirow{ 2}{*}{ 3000} &       iARC &      34 &&          34 &&       324 &      296    \\ 
&&  iR\_Newton &      34 &(\%100) &          34 &(\%100) &       162 &        0    \\ \hline 
\multirow{ 2}{*}{  DIXMAANK} &  \multirow{ 2}{*}{ 3000} &       iARC &      30 &&          30 &&       248 &      225    \\ 
&&  iR\_Newton &      30 &(\%100) &          30 &(\%100) &       124 &        0    \\ \hline 
\multirow{ 2}{*}{  DIXMAANL} &  \multirow{ 2}{*}{ 3000} &       iARC &      29 &&          29 &&       180 &      148    \\ 
&&  iR\_Newton &      29 &(\%100) &          29 &(\%100) &        90 &        0    \\ \hline 
\multirow{ 2}{*}{  DIXMAANM} &  \multirow{ 2}{*}{ 3000} &       iARC &     375 &&         375 &&     10902 &    11542    \\ 
&&  iR\_Newton &     398 &(\%100) &         398 &(\%100) &      6126 &        0    \\ \hline 
\multirow{ 2}{*}{  DIXMAANN} &  \multirow{ 2}{*}{ 3000} &       iARC &      82 &&          82 &&      1368 &     1358    \\ 
&&  iR\_Newton &      87 &(\%100) &          87 &(\%100) &       789 &        0    \\ \hline 
\multirow{ 2}{*}{  DIXMAANO} &  \multirow{ 2}{*}{ 3000} &       iARC &      63 &&          63 &&       908 &      893    \\ 
&&  iR\_Newton &      59 &(\%100) &          59 &(\%100) &       371 &        0    \\ \hline 
\multirow{ 2}{*}{  DIXMAANP} &  \multirow{ 2}{*}{ 3000} &       iARC &      51 &&          51 &&       476 &      432    \\ 
&&  iR\_Newton &      51 &(\%100) &          51 &(\%100) &       238 &        0    \\ \hline 
\multirow{ 2}{*}{  DIXON3DQ} &  \multirow{ 2}{*}{10000} &       iARC &    2257 &&        2256 &&    143968 &   164858    \\ 
&&  iR\_Newton &    2476 &(\%100) &        2476 &(\%100) &     81042 &        0    \\ \hline 
\multirow{ 2}{*}{      DJTL} &  \multirow{ 2}{*}{    2} &       iARC &     215 &&         120 &&       642 &      866    \\ 
&&  iR\_Newton &     204 &(\%32) &          81 &(\%5) &       404 &      512    \\ \hline 
\multirow{ 2}{*}{   DQDRTIC} &  \multirow{ 2}{*}{ 5000} &       iARC &       6 &&           6 &&        34 &       32    \\ 
&&  iR\_Newton &       4 &(\%100) &           4 &(\%100) &        10 &        0    \\ \hline 
\multirow{ 2}{*}{    DQRTIC} &  \multirow{ 2}{*}{ 5000} &       iARC &      15 &&          15 &&        30 &       15    \\ 
&&  iR\_Newton &      11 &(\%100) &          11 &(\%100) &        11 &        0    \\ \hline 
\multirow{ 2}{*}{   EDENSCH} &  \multirow{ 2}{*}{ 2000} &       iARC &      15 &&          15 &&        44 &       25    \\ 
&&  iR\_Newton &      15 &(\%100) &          15 &(\%100) &        22 &        0    \\ \hline 
\multirow{ 2}{*}{       EG2} &  \multirow{ 2}{*}{ 1000} &       iARC &       3 &&           3 &&         6 &        3    \\ 
&&  iR\_Newton &       3 &(\%100) &           3 &(\%100) &         3 &        0    \\ \hline 
\multirow{ 2}{*}{  EIGENALS} &  \multirow{ 2}{*}{ 2550} &       iARC &     179 &&         134 &&     15548 &    20388    \\ 
&&  iR\_Newton &     173 &(\%84) &         150 &(\%89) &      7871 &     1999    \\ \hline 
\multirow{ 2}{*}{   ENGVAL1} &  \multirow{ 2}{*}{ 5000} &       iARC &       9 &&           9 &&        64 &       54    \\ 
&&  iR\_Newton &       9 &(\%100) &           9 &(\%100) &        32 &        0    \\ \hline 
\multirow{ 2}{*}{   ENGVAL2} &  \multirow{ 2}{*}{    3} &       iARC &      21 &&          15 &&       100 &      139    \\ 
&&  iR\_Newton &      21 &(\%57) &          15 &(\%80) &        56 &       34    \\ \hline 
\multirow{ 2}{*}{  ERRINROS} &  \multirow{ 2}{*}{   50} &       iARC &     131 &&         122 &&      1202 &     1106    \\ 
&&  iR\_Newton &     108 &(\%94) &         103 &(\%97) &       504 &       37    \\ \hline 
\multirow{ 2}{*}{  ERRINRSM} &  \multirow{ 2}{*}{   50} &       iARC &     404 &&         396 &&      6566 &     7225    \\ 
&&  iR\_Newton &     167 &(\%98) &         163 &(\%99) &      1154 &       27    \\ \hline 
\multirow{ 2}{*}{    EXPFIT} &  \multirow{ 2}{*}{    2} &       iARC &      14 &&           9 &&        42 &       62    \\ 
&&  iR\_Newton &      11 &(\%27) &           6 &(\%50) &        26 &       38    \\ \hline 
\multirow{ 2}{*}{  EXTROSNB} &  \multirow{ 2}{*}{ 1000} &       iARC &     179 &&         107 &&      2978 &     3586    \\ 
&&  iR\_Newton &     185 &(\%62) &         114 &(\%64) &      1576 &     1553    \\ \hline 
\multirow{ 2}{*}{  FLETBV3M} &  \multirow{ 2}{*}{ 5000} &       iARC &      41 &&          34 &&        86 &       43    \\ 
&&  iR\_Newton &      56 &(\%43) &          32 &(\%41) &        65 &       32    \\ \hline 
\multirow{ 2}{*}{  FLETCHCR} &  \multirow{ 2}{*}{ 1000} &       iARC &    2437 &&        1450 &&     66056 &    90373    \\ 
&&  iR\_Newton &    2187 &(\%66) &        1438 &(\%69) &     29012 &    23819    \\ \hline 
\multirow{ 2}{*}{  FMINSRF2} &  \multirow{ 2}{*}{ 5625} &       iARC &     875 &&         567 &&      6528 &     7378    \\ 
&&  iR\_Newton &     905 &(\%50) &         448 &(\%40) &      2666 &     1989    \\ \hline 
\multirow{ 2}{*}{  FREUROTH} &  \multirow{ 2}{*}{ 5000} &       iARC &      17 &&          11 &&       102 &      120    \\ 
&&  iR\_Newton &      18 &(\%39) &          10 &(\%60) &        51 &       35    \\ \hline 
\multirow{ 2}{*}{  GENHUMPS} &  \multirow{ 2}{*}{ 5000} &       iARC &   14931 &&       11710 &&    477824 &   1724919    \\ 
&&  iR\_Newton &    3567 &(\%2) &        2077 &(\%1) &     25952 &    85744    \\ \hline 
\multirow{ 2}{*}{   GENROSE} &  \multirow{ 2}{*}{  500} &       iARC &     593 &&         350 &&     24494 &    54811    \\ 
&&  iR\_Newton &     690 &(\%19) &         341 &(\%4) &     11862 &    30583    \\ \hline 
\multirow{ 2}{*}{  GROWTHLS} &  \multirow{ 2}{*}{    3} &       iARC &       8 &&           8 &&        32 &       32    \\ 
&&  iR\_Newton &       8 &(\%100) &           8 &(\%100) &        16 &        0    \\ \hline 
\multirow{ 2}{*}{      GULF} &  \multirow{ 2}{*}{    3} &       iARC &      46 &&          31 &&       196 &      249    \\ 
&&  iR\_Newton &      40 &(\%62) &          29 &(\%62) &       101 &       78    \\ \hline 
\multirow{ 2}{*}{     HAIRY} &  \multirow{ 2}{*}{    2} &       iARC &      28 &&          15 &&        84 &      133    \\ 
&&  iR\_Newton &      19 &(\%21) &          10 &(\%40) &        45 &       64    \\ \hline 
\multirow{ 2}{*}{   HATFLDD} &  \multirow{ 2}{*}{    3} &       iARC &      23 &&          19 &&       108 &      140    \\ 
&&  iR\_Newton &      23 &(\%65) &          18 &(\%78) &        64 &       33    \\ \hline 
\multirow{ 2}{*}{   HATFLDE} &  \multirow{ 2}{*}{    3} &       iARC &      23 &&          18 &&       112 &      153    \\ 
&&  iR\_Newton &      24 &(\%62) &          18 &(\%83) &        69 &       45    \\ \hline 
\multirow{ 2}{*}{  HATFLDFL} &  \multirow{ 2}{*}{    3} &       iARC &     851 &&         711 &&      4255 &     5134    \\ 
&&  iR\_Newton &    1127 &(\%84) &         961 &(\%85) &      3429 &     1313    \\ \hline 
\multirow{ 2}{*}{  HEART6LS} &  \multirow{ 2}{*}{    6} &       iARC &    1502 &&         895 &&     15002 &    27063    \\ 
&&  iR\_Newton &    1071 &(\%52) &         620 &(\%51) &      5164 &     6637    \\ \hline 
\multirow{ 2}{*}{  HEART8LS} &  \multirow{ 2}{*}{    8} &       iARC &     115 &&          69 &&      1407 &     2400    \\ 
&&  iR\_Newton &     186 &(\%31) &          97 &(\%23) &      1178 &     2185    \\ \hline 
\multirow{ 2}{*}{     HELIX} &  \multirow{ 2}{*}{    3} &       iARC &      11 &&           7 &&        52 &       72    \\ 
&&  iR\_Newton &      15 &(\%33) &           9 &(\%44) &        48 &       57    \\ \hline 
\multirow{ 2}{*}{  HILBERTA} &  \multirow{ 2}{*}{    2} &       iARC &       5 &&           5 &&        12 &        9    \\ 
&&  iR\_Newton &       3 &(\%100) &           3 &(\%100) &         4 &        0    \\ \hline 
\multirow{ 2}{*}{  HILBERTB} &  \multirow{ 2}{*}{   10} &       iARC &       4 &&           4 &&        16 &       12    \\ 
&&  iR\_Newton &       3 &(\%100) &           3 &(\%100) &         6 &        0    \\ \hline 
\multirow{ 2}{*}{  HIMMELBB} &  \multirow{ 2}{*}{    2} &       iARC &      10 &&           6 &&        26 &       36    \\ 
&&  iR\_Newton &      11 &(\%27) &           6 &(\%33) &        20 &       27    \\ \hline 
\multirow{ 2}{*}{  HIMMELBF} &  \multirow{ 2}{*}{    4} &       iARC &      53 &&          36 &&       310 &      408    \\ 
&&  iR\_Newton &      70 &(\%71) &          62 &(\%77) &       244 &      185    \\ \hline 
\multirow{ 2}{*}{  HIMMELBG} &  \multirow{ 2}{*}{    2} &       iARC &       7 &&           6 &&        21 &       25    \\ 
&&  iR\_Newton &       7 &(\%57) &           6 &(\%67) &        13 &        3    \\ \hline 
\multirow{ 2}{*}{  HIMMELBH} &  \multirow{ 2}{*}{    2} &       iARC &       5 &&           4 &&        13 &       12    \\ 
&&  iR\_Newton &       6 &(\%67) &           4 &(\%75) &         7 &        2    \\ \hline 
\multirow{ 2}{*}{     HUMPS} &  \multirow{ 2}{*}{    2} &       iARC &     125 &&          80 &&       370 &      655    \\ 
&&  iR\_Newton &      89 &(\%13) &          49 &(\%10) &       218 &      291    \\ \hline 
\multirow{ 2}{*}{  HYDC20LS} &  \multirow{ 2}{*}{   99} &       iARC &      11 &&           9 &&       402 &      539    \\ 
&&  iR\_Newton &      11 &(\%73) &           9 &(\%89) &       215 &      165    \\ \hline 
\multirow{ 2}{*}{    JENSMP} &  \multirow{ 2}{*}{    2} &       iARC &       8 &&           8 &&        24 &       26    \\ 
&&  iR\_Newton &       8 &(\%100) &           8 &(\%100) &        16 &        0    \\ \hline 
\multirow{ 2}{*}{    JIMACK} &  \multirow{ 2}{*}{ 3549} &       iARC &      54 &&          54 &&     36564 &    45769    \\ 
&&  iR\_Newton &      52 &(\%100) &          52 &(\%100) &     16267 &        0    \\ \hline 
\multirow{ 2}{*}{    KOWOSB} &  \multirow{ 2}{*}{    4} &       iARC &      20 &&          20 &&       114 &      121    \\ 
&&  iR\_Newton &      18 &(\%94) &          18 &(\%94) &        55 &        4    \\ \hline 
\multirow{ 2}{*}{   LIARWHD} &  \multirow{ 2}{*}{ 5000} &       iARC &      12 &&          12 &&        46 &       45    \\ 
&&  iR\_Newton &      11 &(\%100) &          11 &(\%100) &        21 &        0    \\ \hline 
\multirow{ 2}{*}{  LOGHAIRY} &  \multirow{ 2}{*}{    2} &       iARC &     167 &&         116 &&       390 &      414    \\ 
&&  iR\_Newton &     326 &(\%39) &         233 &(\%49) &       542 &      460    \\ \hline 
\multirow{ 2}{*}{   MANCINO} &  \multirow{ 2}{*}{  100} &       iARC &       6 &&           6 &&        14 &        8    \\ 
&&  iR\_Newton &       4 &(\%100) &           4 &(\%100) &         5 &        0    \\ \hline 
\multirow{ 2}{*}{  MARATOSB} &  \multirow{ 2}{*}{    2} &       iARC &       3 &&           3 &&         7 &        5    \\ 
&&  iR\_Newton &       3 &(\%100) &           3 &(\%100) &         4 &        0    \\ \hline 
\multirow{ 2}{*}{    MEXHAT} &  \multirow{ 2}{*}{    2} &       iARC &      11 &&          11 &&        33 &       43    \\ 
&&  iR\_Newton &      11 &(\%100) &          11 &(\%100) &        22 &        0    \\ \hline 
\multirow{ 2}{*}{    MEYER3} &  \multirow{ 2}{*}{    3} &       iARC &      12 &&          12 &&        48 &       51    \\ 
&&  iR\_Newton &      16 &(\%75) &          15 &(\%73) &        38 &       20    \\ \hline 
\multirow{ 2}{*}{  MODBEALE} &  \multirow{ 2}{*}{20000} &       iARC & ---        &&   --- &&   --- &   ---    \\ 
&&  iR\_Newton &    3317 &(\%99) &        3304 &(\%100) &     65293 &      351    \\ \hline 
\multirow{ 2}{*}{    MOREBV} &  \multirow{ 2}{*}{ 5000} &       iARC &       4 &&           4 &&      1102 &     2064    \\ 
&&  iR\_Newton &       1 &(\%100) &           1 &(\%100) &       401 &        0    \\ \hline 
\multirow{ 2}{*}{  MSQRTALS} &  \multirow{ 2}{*}{ 1024} &       iARC &      39 &&          33 &&      9830 &    12602    \\ 
&&  iR\_Newton &      44 &(\%73) &          36 &(\%83) &      4743 &      149    \\ \hline 
\multirow{ 2}{*}{  MSQRTBLS} &  \multirow{ 2}{*}{ 1024} &       iARC &      32 &&          26 &&      5822 &     7090    \\ 
&&  iR\_Newton &      39 &(\%69) &          31 &(\%81) &      3131 &      156    \\ \hline 
\multirow{ 2}{*}{     NCB20} &  \multirow{ 2}{*}{ 5010} &       iARC &     106 &&          70 &&      2888 &     4664    \\ 
&&  iR\_Newton &      65 &(\%32) &          43 &(\%42) &       688 &      614    \\ \hline 
\multirow{ 2}{*}{    NCB20B} &  \multirow{ 2}{*}{ 5000} &       iARC &      29 &&          18 &&      4286 &     9958    \\ 
&&  iR\_Newton &      38 &(\%47) &          19 &(\%42) &      3297 &     8386    \\ \hline 
\multirow{ 2}{*}{  NONCVXU2} &  \multirow{ 2}{*}{ 5000} &       iARC &   10302 &&       10302 &&     20604 &    10302    \\ 
&&  iR\_Newton &   11094 &(\%100) &       11094 &(\%100) &     11094 &        0    \\ \hline 
\multirow{ 2}{*}{  NONCVXUN} &  \multirow{ 2}{*}{ 5000} &       iARC &   23771 &&       23771 &&     47542 &    23771    \\ 
&&  iR\_Newton &   20913 &(\%100) &       20913 &(\%100) &     20913 &        0    \\ \hline 
\multirow{ 2}{*}{    NONDIA} &  \multirow{ 2}{*}{ 5000} &       iARC &       2 &&           2 &&         4 &        2    \\ 
&&  iR\_Newton &       2 &(\%100) &           2 &(\%100) &         2 &        0    \\ \hline 
\multirow{ 2}{*}{  NONDQUAR} &  \multirow{ 2}{*}{ 5000} &       iARC &      45 &&          37 &&       156 &      126    \\ 
&&  iR\_Newton &      38 &(\%95) &          36 &(\%97) &        70 &        2    \\ \hline 
\multirow{ 2}{*}{  OSBORNEA} &  \multirow{ 2}{*}{    5} &       iARC &      36 &&          28 &&       289 &      412    \\ 
&&  iR\_Newton &      21 &(\%43) &          12 &(\%67) &        75 &       73    \\ \hline 
\multirow{ 2}{*}{  OSBORNEB} &  \multirow{ 2}{*}{   11} &       iARC &      25 &&          19 &&       396 &      539    \\ 
&&  iR\_Newton &      28 &(\%75) &          23 &(\%78) &       225 &      105    \\ \hline 
\multirow{ 2}{*}{  OSCIGRAD} &  \multirow{ 2}{*}{100000} &       iARC &      13 &&          10 &&       190 &      220    \\ 
&&  iR\_Newton &      15 &(\%40) &           9 &(\%56) &        92 &       61    \\ \hline 
\multirow{ 2}{*}{  OSCIPATH} &  \multirow{ 2}{*}{   10} &       iARC &  222974 &&      131997 &&   4233806 &   5617938    \\ 
&&  iR\_Newton &  227426 &(\%59) &      134306 &(\%59) &   2273151 &   2521354    \\ \hline 
\multirow{ 2}{*}{  PALMER1C} &  \multirow{ 2}{*}{    8} &       iARC &     161 &&         161 &&       482 &      362    \\ 
&&  iR\_Newton &      74 &(\%100) &          74 &(\%100) &       145 &        0    \\ \hline 
\multirow{ 2}{*}{  PALMER1D} &  \multirow{ 2}{*}{    7} &       iARC &    1069 &&        1069 &&      3586 &     2567    \\ 
&&  iR\_Newton &     196 &(\%100) &         196 &(\%100) &       379 &        0    \\ \hline 
\multirow{ 2}{*}{  PALMER2C} &  \multirow{ 2}{*}{    8} &       iARC &     109 &&         109 &&       326 &      245    \\ 
&&  iR\_Newton &      76 &(\%100) &          76 &(\%100) &       147 &        0    \\ \hline 
\multirow{ 2}{*}{  PALMER3C} &  \multirow{ 2}{*}{    8} &       iARC &      64 &&          64 &&       252 &      201    \\ 
&&  iR\_Newton &      36 &(\%100) &          36 &(\%100) &        69 &        0    \\ \hline 
\multirow{ 2}{*}{  PALMER4C} &  \multirow{ 2}{*}{    8} &       iARC &      27 &&          27 &&       102 &       84    \\ 
&&  iR\_Newton &      90 &(\%100) &          90 &(\%100) &       177 &        0    \\ \hline 
\multirow{ 2}{*}{  PALMER5C} &  \multirow{ 2}{*}{    6} &       iARC &       9 &&           9 &&        22 &       13    \\ 
&&  iR\_Newton &      10 &(\%100) &          10 &(\%100) &        13 &        0    \\ \hline 
\multirow{ 2}{*}{  PALMER6C} &  \multirow{ 2}{*}{    8} &       iARC &     238 &&         238 &&       870 &      654    \\ 
&&  iR\_Newton &     252 &(\%100) &         252 &(\%100) &       503 &        0    \\ \hline 
\multirow{ 2}{*}{  PALMER7C} &  \multirow{ 2}{*}{    8} &       iARC &      65 &&          65 &&       196 &      143    \\ 
&&  iR\_Newton &      65 &(\%100) &          65 &(\%100) &       120 &        0    \\ \hline 
\multirow{ 2}{*}{  PALMER8C} &  \multirow{ 2}{*}{    8} &       iARC &      76 &&          76 &&       300 &      229    \\ 
&&  iR\_Newton &      90 &(\%100) &          90 &(\%100) &       174 &        0    \\ \hline 
\multirow{ 2}{*}{    PARKCH} &  \multirow{ 2}{*}{   15} &       iARC &      31 &&          22 &&       478 &      685    \\ 
&&  iR\_Newton &      33 &(\%61) &          23 &(\%74) &       270 &      209    \\ \hline 
\multirow{ 2}{*}{  PENALTY1} &  \multirow{ 2}{*}{ 1000} &       iARC &      14 &&          14 &&        28 &       14    \\ 
&&  iR\_Newton &      12 &(\%100) &          12 &(\%100) &        12 &        0    \\ \hline 
\multirow{ 2}{*}{  PENALTY2} &  \multirow{ 2}{*}{  200} &       iARC &      22 &&          22 &&       314 &      319    \\ 
&&  iR\_Newton &      22 &(\%100) &          22 &(\%100) &       157 &        0    \\ \hline 
\multirow{ 2}{*}{  PENALTY3} &  \multirow{ 2}{*}{  200} &       iARC &      24 &&          20 &&       168 &      161    \\ 
&&  iR\_Newton &      25 &(\%60) &          18 &(\%72) &        90 &       44    \\ \hline 
\multirow{ 2}{*}{  POWELLSG} &  \multirow{ 2}{*}{ 5000} &       iARC &      17 &&          17 &&        98 &       91    \\ 
&&  iR\_Newton &      17 &(\%100) &          17 &(\%100) &        49 &        0    \\ \hline 
\multirow{ 2}{*}{    QUARTC} &  \multirow{ 2}{*}{ 5000} &       iARC &      15 &&          15 &&        30 &       15    \\ 
&&  iR\_Newton &      11 &(\%100) &          11 &(\%100) &        11 &        0    \\ \hline 
\multirow{ 2}{*}{   ROSENBR} &  \multirow{ 2}{*}{    2} &       iARC &      29 &&          20 &&        87 &      116    \\ 
&&  iR\_Newton &      32 &(\%62) &          20 &(\%70) &        64 &       51    \\ \hline 
\multirow{ 2}{*}{      S308} &  \multirow{ 2}{*}{    2} &       iARC &      12 &&           9 &&        36 &       41    \\ 
&&  iR\_Newton &      10 &(\%80) &           8 &(\%88) &        20 &        8    \\ \hline 
\multirow{ 2}{*}{  SCHMVETT} &  \multirow{ 2}{*}{ 5000} &       iARC &       5 &&           5 &&       142 &      166    \\ 
&&  iR\_Newton &       6 &(\%100) &           6 &(\%100) &        89 &        0    \\ \hline 
\multirow{ 2}{*}{   SENSORS} &  \multirow{ 2}{*}{  100} &       iARC &      17 &&          12 &&       146 &      263    \\ 
&&  iR\_Newton &      21 &(\%19) &          12 &(\%33) &        66 &       85    \\ \hline 
\multirow{ 2}{*}{   SINEVAL} &  \multirow{ 2}{*}{    2} &       iARC &      66 &&          42 &&       194 &      256    \\ 
&&  iR\_Newton &      63 &(\%63) &          41 &(\%68) &       123 &       81    \\ \hline 
\multirow{ 2}{*}{   SINQUAD} &  \multirow{ 2}{*}{ 5000} &       iARC &      16 &&          11 &&        64 &       63    \\ 
&&  iR\_Newton &      15 &(\%33) &           9 &(\%44) &        32 &       26    \\ \hline 
\multirow{ 2}{*}{    SISSER} &  \multirow{ 2}{*}{    2} &       iARC &      12 &&          12 &&        24 &       12    \\ 
&&  iR\_Newton &      12 &(\%100) &          12 &(\%100) &        12 &        0    \\ \hline 
\multirow{ 2}{*}{     SNAIL} &  \multirow{ 2}{*}{    2} &       iARC &     103 &&          63 &&       290 &      364    \\ 
&&  iR\_Newton &     107 &(\%55) &          63 &(\%56) &       203 &      182    \\ \hline 
\multirow{ 2}{*}{  SPARSINE} &  \multirow{ 2}{*}{ 5000} &       iARC &     153 &&         143 &&     15246 &    18485    \\ 
&&  iR\_Newton &     188 &(\%88) &         174 &(\%94) &     10745 &      183    \\ \hline 
\multirow{ 2}{*}{  SPARSQUR} &  \multirow{ 2}{*}{10000} &       iARC &      15 &&          15 &&        64 &       49    \\ 
&&  iR\_Newton &      15 &(\%100) &          15 &(\%100) &        32 &        0    \\ \hline 
\multirow{ 2}{*}{  SPMSRTLS} &  \multirow{ 2}{*}{ 4999} &       iARC &      17 &&          15 &&       582 &      761    \\ 
&&  iR\_Newton &      17 &(\%76) &          15 &(\%87) &       275 &        4    \\ \hline 
\multirow{ 2}{*}{  SROSENBR} &  \multirow{ 2}{*}{ 5000} &       iARC &       9 &&           7 &&        36 &       35    \\ 
&&  iR\_Newton &      10 &(\%70) &           7 &(\%86) &        20 &       12    \\ \hline 
\multirow{ 2}{*}{  SSBRYBND} &  \multirow{ 2}{*}{ 5000} &       iARC &      75 &&          45 &&     77075 &   177454    \\ 
&&  iR\_Newton &      39 &(\%38) &          23 &(\%52) &     22010 &    11269    \\ \hline 
\multirow{ 2}{*}{   STRATEC} &  \multirow{ 2}{*}{   10} &       iARC &      74 &&          65 &&       886 &     1069    \\ 
&&  iR\_Newton &      67 &(\%90) &          61 &(\%95) &       413 &       87    \\ \hline 
\multirow{ 2}{*}{  TESTQUAD} &  \multirow{ 2}{*}{ 5000} &       iARC &     162 &&         162 &&     16908 &    19812    \\ 
&&  iR\_Newton &     163 &(\%100) &         163 &(\%100) &      8271 &        0    \\ \hline 
\multirow{ 2}{*}{  TOINTGOR} &  \multirow{ 2}{*}{   50} &       iARC &      11 &&          11 &&       234 &      275    \\ 
&&  iR\_Newton &      11 &(\%100) &          11 &(\%100) &       117 &        0    \\ \hline 
\multirow{ 2}{*}{  TOINTGSS} &  \multirow{ 2}{*}{ 5000} &       iARC &       4 &&           4 &&        14 &       10    \\ 
&&  iR\_Newton &       3 &(\%100) &           3 &(\%100) &         7 &        0    \\ \hline 
\multirow{ 2}{*}{  TOINTPSP} &  \multirow{ 2}{*}{   50} &       iARC &      35 &&          22 &&       254 &      335    \\ 
&&  iR\_Newton &      41 &(\%49) &          20 &(\%40) &       156 &       66    \\ \hline 
\multirow{ 2}{*}{  TOINTQOR} &  \multirow{ 2}{*}{   50} &       iARC &       7 &&           7 &&       104 &      103    \\ 
&&  iR\_Newton &       7 &(\%100) &           7 &(\%100) &        52 &        0    \\ \hline 
\multirow{ 2}{*}{  TQUARTIC} &  \multirow{ 2}{*}{ 5000} &       iARC &      11 &&          11 &&        44 &       50    \\ 
&&  iR\_Newton &       1 &(\%100) &           1 &(\%100) &         2 &        0    \\ \hline 
\multirow{ 2}{*}{    TRIDIA} &  \multirow{ 2}{*}{ 5000} &       iARC &      16 &&          16 &&      2128 &     2630    \\ 
&&  iR\_Newton &      17 &(\%100) &          17 &(\%100) &      1310 &        0    \\ \hline 
\multirow{ 2}{*}{    VARDIM} &  \multirow{ 2}{*}{  200} &       iARC &      12 &&          12 &&        24 &       12    \\ 
&&  iR\_Newton &      12 &(\%100) &          12 &(\%100) &        12 &        0    \\ \hline 
\multirow{ 2}{*}{  VAREIGVL} &  \multirow{ 2}{*}{   50} &       iARC &       5 &&           5 &&        42 &       38    \\ 
&&  iR\_Newton &       5 &(\%100) &           5 &(\%100) &        21 &        0    \\ \hline 
\multirow{ 2}{*}{  VIBRBEAM} &  \multirow{ 2}{*}{    8} &       iARC &      70 &&          41 &&       644 &     1131    \\ 
&&  iR\_Newton &      39 &(\%31) &          25 &(\%48) &       183 &      226    \\ \hline 
\multirow{ 2}{*}{    WATSON} &  \multirow{ 2}{*}{   12} &       iARC &      14 &&          14 &&       174 &      214    \\ 
&&  iR\_Newton &      14 &(\%100) &          14 &(\%100) &        88 &        0    \\ \hline 
\multirow{ 2}{*}{     WOODS} &  \multirow{ 2}{*}{ 4000} &       iARC &      15 &&          15 &&        40 &       26    \\ 
&&  iR\_Newton &     172 &(\%87) &         157 &(\%92) &       404 &      144    \\ \hline 
\multirow{ 2}{*}{     YFITU} &  \multirow{ 2}{*}{    3} &       iARC &      54 &&          38 &&       270 &      348    \\ 
&&  iR\_Newton &      55 &(\%65) &          39 &(\%64) &       172 &      130    \\ \hline 
\multirow{ 2}{*}{  ZANGWIL2} &  \multirow{ 2}{*}{    2} &       iARC &       3 &&           3 &&         6 &        3    \\ 
&&  iR\_Newton &       1 &(\%100) &           1 &(\%100) &         1 &        0    \\ \hline 
\end{longtable}
} 
\end{scriptsize}


\bibliographystyle{IMANUM-BIB}
\bibliography{itrace_references}

\end{document}